\title{\bf Very large graphs}
\author{\sc L\'aszl\'o Lov\'asz\footnote{Research sponsored by OTKA Grant No.~67867.}}
\date{December 2008\\
\begin{flushright}
\it Dedicated to the memory of Oded Schramm\\[1cm]
\end{flushright}
}

\documentclass[fleqn]{article}
\usepackage{graphicx,amssymb,amsmath,theorem,hyperref}

\long\def\killtext#1{}
\newtheorem{theorem}{Theorem}[section]
\newtheorem{prop}[theorem]{Proposition}
\newtheorem{lemma}[theorem]{Lemma}

\newtheorem{corollary}[theorem]{Corollary}

\theorembodyfont{\rmfamily}
\newtheorem{remark}[theorem]{Remark}
\newtheorem{example}[theorem]{Example}
\newtheorem{alg}[theorem]{Algorithm}
\newtheorem{prob}[theorem]{Problem}
\newtheorem{conj}[theorem]{Conjecture}
\newenvironment{proof}{\noindent{\bf Proof. }}{\hfill$\square$\medskip}
\begin{document}

\addtolength{\baselineskip}{3pt}

\def\R{{\mathbb R}}
\def\one{{\mathbf 1}}
\def\Q{{\mathbb Q}}
\def\Z{{\mathbb Z}}
\def\N{{\mathbb N}}
\def\C{{\mathbb C}}
\def\hom{{\rm hom}}
\def\Hom{{\rm Hom}}
\def\Inj{{\rm Inj}}
\def\Ind{{\rm Ind}}
\def\inj{{\rm inj}}
\def\ind{{\rm ind}}
\def\sur{{\rm sur}}
\def\PAG{{\rm PAG}}
\def\eps{\varepsilon}
\def\Ge{{\mathbb G}}
\def\Gb{{\mathbf G}}
\def\IO{{\infty\to1}}
\def\GR{\text{\rm GR}}
\def\irreg{\text{\rm irreg}}
\def\Mch{\text{\sf match}}
\def\col{\text{\rm col}}
\def\maxcut{\text{\sf maxcut}}
\def\Maxcut{\text{\sf Maxcut}}
\def\mmcut{\text{\sf mmcut}}
\def\rmcut{\text{\sf rmcut}}

\def\proofend{\hfill$\square$\medskip}
\def\Proof{\noindent{\bf Proof. }}

\def\pm{\text{\sf pm}}

\def\AA{{\cal A}}\def\BB{{\cal B}}\def\CC{{\cal C}}\def\DD{{\cal D}}\def\EE{{\cal E}}\def\FF{{\cal F}}
\def\GG{{\cal G}}\def\HH{{\cal H}}\def\II{{\cal I}}\def\JJ{{\cal J}}\def\KK{{\cal K}}\def\LL{{\cal L}}
\def\MM{{\cal M}}\def\NN{{\cal N}}\def\OO{{\cal O}}\def\PP{{\cal P}}\def\QQ{{\cal Q}}\def\RR{{\cal R}}
\def\SS{{\cal S}}\def\TT{{\cal T}}\def\UU{{\cal U}}\def\VV{{\cal V}}\def\WW{{\cal W}}\def\XX{{\cal X}}
\def\YY{{\cal Y}}\def\ZZ{{\cal Z}}
\def\QG{{\cal QG}}
\def\QGM{{\cal QGM}}
\def\CD{{\cal CD}}
\def\Prob{{\sf P}}
\def\E{{\sf E}}
\def\Var{{\sf Var}}
\def\T{^{\sf T}}

\def\tr{{\rm tr}}
\def\cost{\hbox{\rm cost}}
\def\val{\hbox{\rm val}}
\def\rk{{\rm rk}}
\def\diam{{\rm diam}}
\def\simi{{\rm sim}}
\def\diag{{\rm diag}}
\def\Ker{{\rm Ker}}

\maketitle

\tableofcontents

\section{Introduction}

\subsection{Huge networks}

In the last decade it became apparent that a large number of the most
interesting structures and phenomena of the world can be described by
networks: separable elements, with connections (or interactions)
between certain pairs of them.

\begin{itemize}
\item Among such a networks, the best known and the most studied is
the {\it internet}. Moreover, the internet (as the physical
underlying network) gives rise to many of the networks: the network
of hyperlinks (web, logical Internet), Internet based social
networks, distributed data bases, etc. The size of the internet is
growing fast: currently the number of web pages may be 30 billion or
more, and the number of devices is probably more than a billion.

\item Social networks are basic objects of many studies in the area
of sociology, history, epidemiology and economics. The largest social
network is the acquaintance graph of all living people, with about 7
billion nodes.

\item Biology contributes {\it ecological networks}, networks of
interactions between proteins, and the human brain, just to mention a
few. The human brain is really large for its mass, having about
$10^{11}$ nodes.

\item Statistical physics studies the interactions between large
numbers of discrete particles, where the underlying structure is
often described by a graph. For example, a crystal can be though of
as a graph whose nodes are the atoms and whose edges represent
chemical bonds. A perfect crystal is a rather boring graph, but
impurities and imperfections create interesting graph-theoretical
digressions. 12 gram of a diamond has about $6\times 10^{23}$ nodes.

\item Some of the largest networks in engineering occur in chip
design. Even though these networks are man-made and planned, many of
their properties are difficult to determine by computation due to
their huge size. There can be more than a billion transistors on a
chip now.

\item To be pretentious, we can say that the whole universe is a
single (really huge, possibly infinite) network, where the nodes are
events (interactions between elementary particles), and the edges are
the particles themselves. This is a network with perhaps $10^{80}$
nodes.
\end{itemize}

These huge networks pose exciting challenges for the mathematician.
Graph Theory (the mathematical theory of networks) has been one of
the fastest developing areas of mathematics in the last decades; with
the appearance of the Internet, however, it faces fairly novel,
unconventional problems. In traditional graph theoretical problems
the whole graph is exactly given, and we are looking for
relationships between its parameters or efficient algorithms for
computing its parameters. On the other hand, very large networks
(like the Internet) are never completely known, in most cases they
are not even well defined. Data about them can be collected only by
indirect means like random local sampling or by monitoring the
behavior of various global processes.

Dense networks (in which a node is adjacent to a positive percent of
others nodes) and sparse networks (in which a node has a bounded
number of neighbors) show a very diverse behavior. From a practical
point of view, sparse networks are more important, but at present we
have more complete theoretical results for dense networks.

\subsection{What to ask about them?}\label{QUEST}

Let us discuss three possible questions that can be asked about a
really large graph, say the internet.

\medskip

{\bf Question 1.} {\it Does the graph have an odd number of nodes?}

\smallskip

This is a very basic property of a graph in the classical setting.
For example, it is one of the first theorems or exercises in a graph
theory course that every graph with an odd number of nodes has a node
with even degree.

But for the internet, this question is clearly nonsense. Not only
does the number of nodes change all the time, with devices going
online and offline, but even if we fix a specific time like 12:00am
today, it is not well-defined: there will be computers just in the
process of booting up, breaking down etc.

\medskip

{\bf Question 2.} {\it What is the average degree of nodes?}

\smallskip

This, on the other hand, is a meaningful question. Of course, the
average degree can only be determined with a certain error, and it
will change with technology or the social composition of users; but
at a given time, a good approximation can be sought (I am not
speaking now about how to find it).

\medskip

{\bf Question 3.} {\it Is the graph connected?}

\smallskip

To this question, the answer is almost certainly no: somewhere there
will be a faulty router with some unhappy users on the wrong side of
it. But this is not the interesting way to ask the question: we
should consider the internet disconnected if, say, an earthquake
combined with a sunflare severs all connections between the Old and
New worlds. So we want to ignore small components that are negligible
with respect to the whole graph, and consider the graph disconnected
only if it decomposes into two parts which are commeasurable with the
whole. On the other hand, we may want to allow that the two parts be
connected by a few edges, and still consider the graph disconnected.

\medskip

{\bf Question 4.} {\it Find the largest cut in the graph.}

\smallskip

(This means to find the partition of the nodes into two classes so as
to maximize the number of edges connecting the two classes.) This
example shows that even if the question is meaningful, it is not
clear in what form can we expect the answer. The fraction of edges
contained in the largest cut can be determined relatively easily
(with and error that is small with large probability); but how to
specify the largest cut itself (or even an approximate version of
it)?

\subsection{How to obtain information about them?}

If we face a large network (think of the internet) the first
challenge is to obtain information about it. Often, we don't even
know the number of nodes.

\subsubsection{Local sampling}\label{LOC-SAMP}

Properties of very large graphs can be studied by sampling small
subgraphs. The theory of this, called {\it property testing} in
computer science, emerged in the last decade, and will be one of the
main concerns of this paper.

In the case of dense graphs $G$, the sampling process is simple: we
select independently a fixed number $k$ of random nodes, and
determine the edges between them, to get a random induced subgraph.
We'll call this {\it subgraph sampling}. For each graph $F$, this
defines a probability of seeing $F$ when $|V(F)|$ nodes are sampled,
and so it gives a probability distribution $\sigma_{G,k}$ on all
graphs with $k$ (labeled) nodes. It turns out that this sample
contains enough information to determine many properties and
parameters of the graph (with an error that is with large probability
arbitrarily small if $k$ is sufficiently large depending only on the
error bound).

To get a mathematically exact description of algorithms for very
large graphs, we define a {\it subgraph sampling oracle} as a black
box that, for a given positive integer $m$, returns a random $m$-node
graph from some (otherwise unknown) distribution. We think of this as
a random induced subgraph of a very large, otherwise unknown graph
$G$. We assume that the oracle is consistent in the sense that for
any $k$ there is a graph $G$ such that the distribution of the
$k$-samples from $G$ is arbitrarily close to the distribution of the
answers by the oracle. (Theorem \ref{THM:GRAPHON-CHAR} will give a
characterization of consistent distributions.)

In the case of sparse graphs with bounded degree, the subgraph
sampling method gives a trivial result: the sampled subgraph will
almost certainly be edgeless. Probably the most natural way to fix
this is to consider {\it neighborhood sampling}. Let $\GG_d$ denote
the class of finite graphs with all degrees bounded by $d$. For
$G\in\GG_d$, select a random node and explore its neighborhood to a
given depth $m$. This provides a probability distribution
$\rho_{G,m}$ on graphs in $\GG_d$, with a specified root node, such
that all nodes are at distance at most $m$ from the root. We will
shortly refer to these rooted graphs as $m$-balls. Note that the
number of possible $m$-balls is finite if $d$ and $m$ are fixed. We
can formulate this abstractly as a {\it neighborhood sampling
oracle}, a black box that, for a given positive integer $m$, returns
an $m$-ball.

The situation for sparse graphs is, however, less satisfactory than
for dense graphs, for two reasons. First, a full characterization of
consistent neighborhood sampling oracles is not known (cf. Conjecture
\ref{ALD-LYO}). Second, neighborhood sampling does not reveal
important global properties of the graph like expansion. This
suggests looking at further possibilities. Suppose, for example, that
instead of exploring the neighborhood of a single random node, we
could select two (or more) random nodes and determine simple
quantities associated with them, like pairwise distances, maximum
flow, electrical resistance, hitting times of random walks. What
information can be gained by such tests? Is there a ``complete'' set
of tests that would give enough information to determine the global
structure of the graph to a reasonable accuracy? These methods should
lead to different theories of large graphs and their limit objects,
largely unexplored.

Sample distribution (in both the dense and sparse cases) are
equivalent to counting induces subgraphs of a given type. Instead of
this, we could count homomorphism (or injective homomorphisms) of a
``small'' graphs into the graph. The connection with sample
distribution can be expressed by inclusion-exclusion formulas, and it
is not essential. Often homomorphism numbers are algebraically better
behaved, and they also have the advantage that they suggest
different, ``dual'' approaches by reversing the arrows in the
category of graph homomorphisms.

\subsubsection{Observing global processes}

Another source of information about a network is the observation of
the behavior of various global processes either globally (through
measuring some global parameter), or locally (at one node, or a few
neighboring nodes, but for a longer time). Statistical physical
models on the graph are examples of the first kind of approach (we
return to them in section \ref{STATPHYS}). Crawlers can be considered
as examples of the second, and there are some sporadic results about
the local observation of simpler, random processes \cite{BeLo,BKLRT}.
A general theory of such local observation has not emerged yet
though.

\subsubsection{Left and right homomorphisms}

Instead of testing, it is often more convenient to talk about
homomorphisms (adjacency-preserving maps) between graphs. This leads
to the following setup. If we are given a (large) graph $G$, we may
try to study its local structure by counting homomorphisms from
various ``small'' graphs $F$ into $G$; and we can study its global
structure by counting its homomorphisms into various small graphs
$H$. The first type of information is closely related (in many cases,
equivalent) to sampling, while the second is related to statistical
physics. As in statistical physics, one needs weighted graphs $H$
here to get meaningful results.

\subsection{How to model them?}

\subsubsection{Random graphs}

We are celebrating the 50-th birthday of random graphs this year: The
simplest random graph model was developed by Erd\H{o}s and R\'enyi
\cite{ER1} and Gilbert \cite{Gil} in 1959. Given a positive integer
$n$ and a real number $0\le p\le 1$, we generate a random graph
$\Ge(n,p)$ by taking $n$ nodes, say $[n]=\{1,\dots,n\}$, and
connecting any two of them with probability $p$, making an
independent decision about every pair.

There are alternate models, essentially equivalent: we could fix the
number of edges $m$, and then choose a random $m$-element subset of
the set of pairs in $[n]$, uniformly from all such subsets. This
random graph $\Ge(n,m)$ is very similar to $\Ge(n,p)$ when $m=p
\binom{n}{2}$. Another model, closer to some of the more recent
developments, is {\it evolving random graphs}, where edges are added
one by one, always choosing uniformly from the set of unconnected
pairs. Stopping this process after $m$ steps, we get $\Ge(n,m)$.

Erd\H{o}s--R\'enyi random graphs have many interesting, often
surprising properties, and a huge literature, see \cite{Bol1,JLR}.
One conventional wisdom about random graphs with a given edge density
is that they are all alike. Their basic parameters, like chromatic
number, maximum clique, triangle density, spectra etc. are highly
concentrated. This fact will be an important motivation when defining
the right measure of global similarity of graphs.

Many generalizations of this random graph model have been studied.
For example, one could have different probabilities assigned to
different edges. A variation of this idea, discovered independently
in \cite{LSz1}, \cite{BJR} and perhaps elsewhere, is the notion of
$W$-random graphs, to be discussed in section \ref{WRAND} and used
throughout these notes.

\subsubsection{Randomly growing graphs}

Random graph models on a fixed set of nodes, discussed above, fail to
reproduce important properties of real-life networks. For example,
the degrees of Erd\H{o}s--R\'enyi random graphs follow a binomial
distribution, and so they are asymptotically normal if the edge
probability $p$ is a constant, and asymptotically Poisson if the
expected degree is constant (i.e., $p=p(n)\sim c/n$). In either case,
the degrees are highly concentrated around the mean, while the
degrees of real life networks tend to obey the ``Zipf phenomenon'',
which means that the tail of the distribution decreases according to
a power law.

In 2002 Albert and Barab\'asi \cite{AlBa,Bar} created a random
network model growing according to natural rules, which could
reproduce this behavior. Since then a lot of variations of growing
networks were introduced. The process of graph generation usually
consists of random steps obeying some local rules.

This is perhaps the first point which suggests one of our main tools,
namely assigning limits to sequences of graphs. Just as the Law of
Large Numbers tells us that adding up more and more independent
random variables we get an increasingly deterministically behaving
number, these growing graph sequences tend to have a well-defined
structure, independent of the random choices made along the way. In
the limit, the randomness disappears, and the asymptotic behavior of
the sequence can be described by a well-defined limit object. You
will find more on this in Sections \ref{INTRO:LIMITS} and
\ref{LIM-DENSE}.

\subsubsection{Quasirandom graphs}\label{QUASIRAND}

The theory of quasirandom graphs, introduced by Thomason \cite{Tho}
and Chung, Graham and Wilson \cite{CGW}, is based on the following
observation: not only have random graphs a variety of quite strict
properties (with large probability), but for several of these basic
properties, the exceptional graphs are the same. In other words, any
of these properties implies the others, regardless of any stochastic
consideration.

To make this idea precise, we consider a sequence of graphs $(G_n)$
with $|V(G_n)|\to\infty$. For simplicity, assume that $|V(G_n)|=n$.
Let $0<p<1$ be a real number. Consider the following properties of
these graphs.

\smallskip

(P1) All degrees are asymptotically $p n$ and all codegrees (numbers
of common neighbors of two nodes) are asymptotically $p^2 n$.

\smallskip

(P2) For every fixed graph $F$, the number of homomorphisms of $F$
into $G_n$ is asymptotically $p^{|E(F)|} n^{|V(F)|}$.

\smallskip

(P3) The number of edges is asymptotically $pn^2/2$ and the number of
$4$-cycles is asymptotically $p^4 n^4/8$.

\smallskip

(P4) The number of edges induced by a set of nodes of size $\alpha n$
is asymptotically $p \alpha^2n^2/2$.

\smallskip

All these properties hold with probability $1$ if $G_n=\Ge(n,p)$.
However, more is true: if a graph sequence satisfies either one of
them, then it satisfies all \cite{CGW}. Such graph sequences are
called {\it quasirandom}. The four properties above are only a
sampler; there are many other random-like properties that are also
equivalent to these \cite{CGW,SiSo1,SiSo2}.

Many interesting deterministic graph sequences are quasirandom. We
mention an important example from number theory:

\begin{example}\label{EXA:PALEY}
{\it Paley graphs.} Let $p_n$ be the $n$-th prime congruent $1$
modulo $4$, and let us define a graph on $\{1,\dots,p_n\}$ by
connecting $i$ and $j$ if and only if $i-j$ is a quadratic residue.
The Paley graphs converge to the function $W\equiv 1/2$.
\end{example}

The theory of convergent graph sequences (Section
\ref{CONV-LIM-DENSE}) can be considered as a rather far-reaching
generalization of quasirandom sequences.

\subsection{How to approximate them?}

We want a compact approximate description of a very large network,
usually in the form a (relatively) small networks or at least a
network with a compact description. To make this mathematically
precise, we need to define what we mean by two graphs to be
``similar'' or ``close'', and describe what kind of structures we use
for approximation.

\subsubsection{The distance of two graphs}

There are many ways of defining the distance of two graphs $G$ and
$G'$. Suppose that the two graphs have a common node set $[n]$. Then
a natural notion of distance is the {\it edit distance}, defined as
the number of edges to be changed to get from one graph to the other.
Since our graphs are very large, we want to normalize this, and
define
\[
d_1(G,G')=\frac{|E(G)\triangle E(G')|}{\binom{n}{2}}.
\]
While this distance plays an important role in the study of testable
graph properties, it does not reflect structural similarity well. To
raise one objection, consider two random graphs on $[n]$ with
edge-density $1/2$. As mentioned in the introduction, these graphs
are very similar from almost every aspect, but their normalized edit
distance is large (about $1/2$ with large probability). One might try
to improve this by relabeling one of them to get the best overly
minimizing the edit distance; but the improvement would be marginal
($o(1)$).

Another trouble with the notion of edit distance is that it is only
defined when the two graphs have the same number of nodes.

We could base the measurement of distance on sampling. We define the
{\it sampling distance} of two graphs $G$ and $G'$ by
\begin{equation}\label{EQ:SAMPDIST}
d_{\rm sample}(G,G')=\sum_{k=1}^{\infty} \frac{1}{2^k} d_{\rm
tv}(\sigma_{G,k},\sigma_{G',k})
\end{equation}
(where $d_{\rm tv}(\alpha,\beta)=\sup_X |\alpha(X)-\beta(X)|$ denotes
the total variation distance of the distributions $\alpha$ and
$\beta$). Here the coefficients $1/2^k$ are quite arbitrary, only to
make the sum convergent. This distance, however, would not directly
reflect any structural similarity.

In section \ref{DIST} we will define a further distance between
graphs, which will be satisfactory from all these points of view: it
will be defined for two graphs with possibly different number of
nodes, the distance of two random graphs with the same edge density
will be very small, and it will reflect global structural similarity.
It will define the same topology as $d_{\rm sample}$.

The construction of the sampling distance can be carried over to
bounded degree graphs, by replacing in \eqref{EQ:SAMPDIST} the
sampling distributions $\sigma_{G,k}$ by the neighborhood
distributions $\rho_{G,k}$. We must point out, however, that it seems
to be difficult to define a notion of distance between two graphs
with bounded degree reflecting global similarity.

\subsubsection{Approximation by smaller: Regularity Lemma}

As the exact description of huge networks is not known, and they are
too big for direct study (e.g., for testing different algorithms or
protocols directly on the whole internet), an important operation
would be to ``scale down'' by producing a smaller network with
similar properties. The main tool for doing so is the
``Szemer\'edi-partition'' or ``regularity Lemma''. Szemer\'edi
developed his Regularity Lemma for his celebrated proof of the
Erd\H{o}s--Tur\'an Conjecture on arithmetic progressions in dense
sets of integers in 1975. Since then, the Lemma has emerged as a
fundamental tool in graph theory, with many applications in extremal
graph theory, combinatorial number theory, graph property testing
etc., and became a true focus of research in the past years.

This lemma can be viewed as an archetypal example of dichotomy
between randomness and structure, where we try to decompose a (large
and complicated) object $A$ into a more highly structured object $A'$
with a (quasi)random perturbation (cf. Tao \cite{Tao}). The highly
structured part may be easier to handle, the quasirandom part will
often be simpler due to Laws of Large Numbers. We'll introduce this
partition in section \ref{SZEMEREDI} (and use it throughout).

Finding the Szemer\'edi partition of a huge dense graph is an example
of the problem posed in Question 4 in Section \ref{QUEST}. Algorithm
\ref{SZEM-ALG} will be an example of a possible solution: how a
partition of the nodes can be determined in an implicit form, even if
describing for each node which class it belongs to would take too
much space.

\subsubsection{Approximation by infinite: convergence and
limits}\label{INTRO:LIMITS}

This idea can be motivated by how we look at a large piece of metal.
This is a crystal, that is a really large graph consisting of atoms
and bonds between them. But from many points of view (e.g., the use
of the metal in building a bridge), it is more useful to consider it
as a continuum with a few important parameters (density, elasticity
etc.). Its behavior is governed by differential equations. Can we
consider a more general very large graph as some kind of continuum?

One way to make this intuition precise is to consider a growing
sequence $(G_n)$ of graphs whose number of nodes tends to infinity,
and to define when such a sequence is convergent. (We have mentioned
this idea in connection with randomly growing graphs, but now we
don't assume anything about how the graphs in the sequence are
obtained.) Our discussion of sampling suggests a general principle
leading to a definition: we consider samples of a fixed size $k$ from
$G_n$, and their distribution. We say that the sequence is {\it
locally convergent} (with respect to the given sampling method) if
this distribution tends to a limit as $n\to\infty$ for every fixed
$k$. The family of limiting distributions (one for each $k$) can be
considered as a limit object of the sequence.

For dense graphs, this notion of convergence was suggested by
Erd\H{o}s, Lov\'asz and Spencer \cite{ELS}, and elaborated by Borgs,
Chayes, Lov\'asz, S\'os, Szegedy and Vesztergombi
\cite{BCLSSV,BCLSV1,BCLSV2}. For sparse graphs, this kind of
convergence was introduced by Aldous \cite{Al} and by Benjamini and
Schramm \cite{BSch}. These notions will be discussed in Sections
\ref{CONV-DENSE} and \ref{CONV-SPARSE}, respectively.

The definition above represents the limit of a graph sequence as a
collection of probability distributions on graphs, one for each
sample size. This is not always a helpful representation of the limit
object, and a more explicit description is desirable. A next step is
to represent the family of distributions on finite graphs (the
samples) by a single probability distribution on countable graphs.
For sparse graphs, Benjamini and Schramm provide such a description
as certain measures on countable rooted graphs with bounded degree
(see section \ref{GRAPHING}, and a similar description for dense
graph limits is also known as certain ergodic measures on countable
graphs (\cite{Sze2}; see Theorem \ref{THM:GRAPHON-CHAR}).

More explicit descriptions of these limit objects can also be given.
Let us start with the dense case. Here the limit object can be
described as a two-variable measurable function $W:~[0,1]^2\to[0,1]$,
called a {\it graphon} (Lov\'asz and Szegedy \cite{LSz1}; see Section
\ref{GRAPHON}). These limit objects can be considered as weighted
graphs with a continuum underlying set, or (if you wish) as graphs on
a nonstandard model of the unit interval.

\begin{figure}[htb]
  \centering
  \includegraphics*[height=120pt]{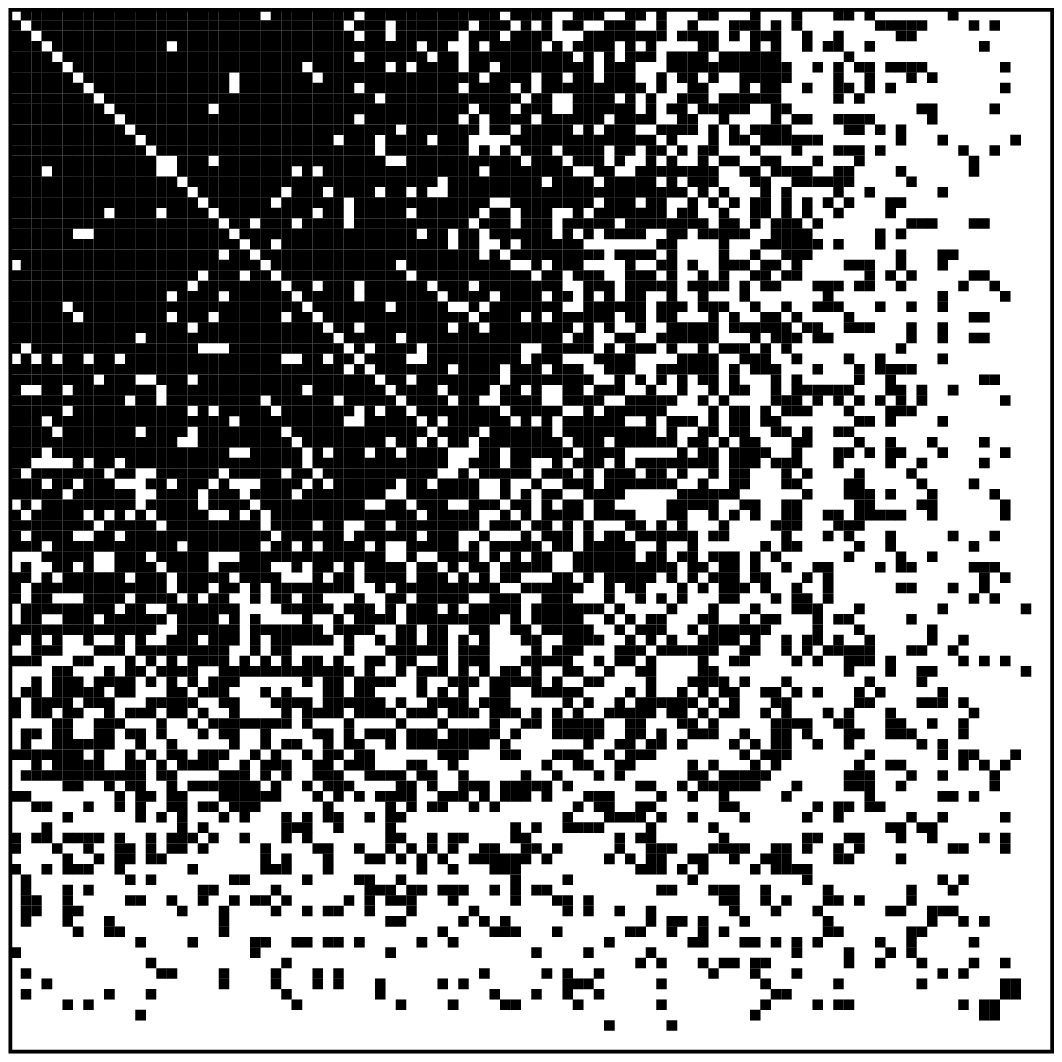}
  \includegraphics*[height=110pt]{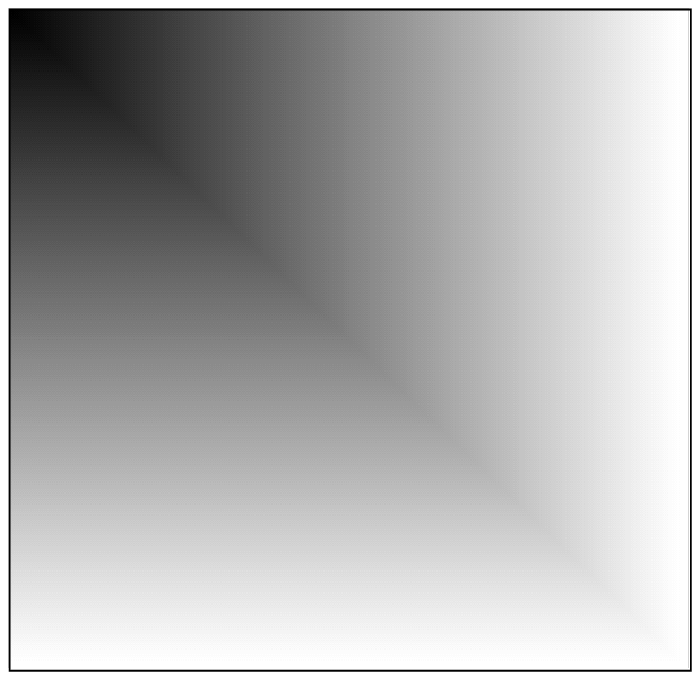}
  \caption{A randomly grown uniform attachment graph with 100 nodes}
  \label{FIG:UNIF}
\end{figure}

Let us describe an example here; more to follow in Section
\ref{LIM-EXAMPLES}. The picture on the left hand side of Figure
\ref{FIG:UNIF} is the adjacency matrix of a graph $G$ with 100 nodes,
where the 1's are represented by black squares and the 0's, by white
squares. The graph itself is constructed by a simple randomized
growing rule: Starting with a single node, we either add a new node
or a new edge; a new node is born with probability $1/n$, where $n$
is the current number of nodes. (A closely related graph sequence
(randomly grown uniform attachment graphs) will be discussed in
detail in Section \ref{LIM-EXAMPLES}.)

The picture on the right hand side is a grayscale image of the
function $U(x,y)=1-\max(x,y)$. The similarity with the picture on the
left is apparent; and suggests that the limit of the graph sequence
on the left is this function. This turns out to be the case in a well
defined sense. It follows that to approximately compute various
parameters of the graph on the left hand side, we can compute related
parameters of the function on the right hand side. For example, the
triangle density of the graph on the left tends (as $n\to\infty$) to
the integral
\[
\int_{[0,1]^3} U(x,y)U(y,z)U(z,x)\,dx\,dy\,dz.
\]

Two more remarks on the dense case. Of course, a graphon can be
infinitely complicated. But in many cases limits of growing graph
sequences have a limit graphon that is a continuous function
described by a simple formula (see a couple of examples in Section
\ref{LIM-EXAMPLES}). Such a limit graphon provides a very useful
approximation of a large dense graph.

Instead of the interval $[0,1]$, we can consider any probability
space $(\Omega,\AA,\pi)$ with a symmetric measurable function
$W:~\Omega\times\Omega\to [0,1]$. This would not give a greater
generality, but it is sometimes useful to represent the limit object
by other probability spaces. We'll see an example of this in Section
\ref{LIM-EXAMPLES}.

In the sparse case, the limit object can be described as a {\it
graphing} (known from group theory or ergodic theory, Elek
\cite{Elek1}), or as a measure preserving graph (see Section
\ref{GRAPHING}), or as a distribution on rooted countable graphs with
special properties.

Instead of sampling, we can use dual (global) measurements, more
precisely, homomorphisms into fixed small graphs, to define
convergence. The remarkable fact is that under the right conditions,
this leads to an equivalent notion! (See Sections \ref{DENSE-CONV-R},
\ref{SP-CONV-R}.)

\subsubsection{Optimization problems for graphs}\label{GR-OPT}

We have presented the theory of convergent graph sequences and their
limits as an answer to problems coming from very large networks, but
a very strong motivation comes from extremal graph theory.

Consider the following two optimization problems.

\smallskip

{\bf Classical optimization problem.} Find the minimum of $x^3-6x$
where $x$ is a nonnegative real number.

\smallskip

{\bf Graph optimization problem.} Find the minimum of $t(C_4,G)$ over
all graphs $G$ with $t(K_2,G)\ge 1/2$. (Here $t(F,G)$, the {\it
homomorphism density} of $F$ in $G$, denotes the probability that a
random map of $V(F)$ into $V(G)$ preserves the edges. $C_4$ denotes
the $4$-cycle and $K_2$ is the complete graph with 2 nodes.)

\smallskip

The solution of the classical optimization problem is of course
$x=\sqrt{2}$. This means that it has no solution over the rationals,
but we can find rational numbers that are arbitrarily close to being
optimal. If we want a single solution, we have to go to the
completion of the rationals, i.e., to the reals.

The graph optimization problem may take a bit more effort to solve,
but it is not hard to show that if the edge-density is $1/2$, then
the $4$-cycle density is larger than $1/16$. Furthermore, this
density gets arbitrarily close to $1/16$ for appropriate families of
graphs: the most important example is a random graph with
edge-density $1/2$ (cf.\ also Section \ref{QUASIRAND} and Theorem
\ref{THM:4CYC}).

This suggests that we could try to enlarge the set of (finite) graphs
with new objects so that the appropriate extension of our
optimization problem has a solution among the new objects.
Furthermore, we want that these new objects should be approximable by
graphs, just like real numbers are approximable by rationals.

\smallskip

Many of the basic tools in the theory of very large graphs have been
first applied in extremal graph theory: the Regularity Lemma
\cite{Szem2}, convergent graph sequences \cite{ELS}, quasirandom
graphs \cite{Tho,CGW}.

The example above shows that limit objects may provide cleaner
formulations of extremal graph theory results, with no error terms.
In some cases this goes further, and the limit objects provide a way
to state, in an exact way, questions like "How do extremal graphs
look like?". They have similar uses in the theory of computing. We
discuss these applications in Sections \ref{TEST} and
\ref{EXT-GRAPHONS}.

\subsection{Mathematical tools}

It is clear from the above that this area is at the crossroads of
different fields of mathematics. Graph theory and computer science
are the main sources, and probability and mathematical statistics are
crucial tools. Group theory, in particular finitely generated groups,
have provided many of the questions and ideas in the theory of limits
of graphs with bounded degree. Ergodic theory may play a similar role
in the dense case. Measure theory is needed, and an important new
general proof method uses nonstandard analysis.

We will discuss one further tool, namely Frobenius algebras, which
are used in the proofs of characterization theorems of homomorphism
functions, but also in some other studies of graph parameters; see
Section \ref{ALGEBRA}.

\section{Graph parameters}\label{GRAPHPARAMS}

A {\it graph parameter} is a real valued function defined on
isomorphism types of graphs (including the graph $K_0$ with no nodes
and edges). A {\it simple graph parameter} is defined only on
isomorphism types of simple graphs (i.e., on graphs with no loops or
multiple edges). A graph parameter $f$ is {\it multiplicative} if
$f(G)=f(G_1)f(G_2)$ whenever $G$ is the disjoint union of $G_1$ and
$G_2$. We say that a graph parameter is {\it normalized} if its value
on $K_1$, the graph with one node and no edge, is $1$. Note that if a
graph parameter is multiplicative and not identically $0$, then its
value on $K_0$ (the graph with no nodes and no edges) is $1$.

\subsection{Connection matrices and reflection positivity}

A {\it $k$-labeled graph} is a graph in which $k$ of the nodes are
labeled by $1,\dots,k$ (there may be any number of unlabeled nodes).
A $0$-labeled graph is just an unlabeled graph.

Let $F_1$ and $F_2$ be two $k$-labeled graphs. We define the
$k$-labeled graph $F_1F_2$ by taking their disjoint union, and then
identifying nodes with the same label. Clearly this multiplication is
associative and commutative. For two $0$-labeled graphs, $F_1F_2$ is
their disjoint union.

Let $f$ be any graph parameter and fix an integer $k\ge 0$. We define
the $k$-th {\it connection matrix} of the graph parameter $f$ as the
(infinite) symmetric matrix $M(f,k)$, whose rows and columns are
indexed by (isomorphism types of) $k$-labeled graphs, and the entry
in the intersection of the row corresponding to $F_1$ and the column
corresponding to $F_2$ is $f(F_1F_2)$.

We call the graph parameter {\it reflection positive} if all the
corresponding connection matrices are positive semidefinite.

\subsection{Homomorphisms from the left}

\subsubsection{Versions of homomorphism numbers}

For two finite graphs $F$ and $G$, let $\hom(F,G)$ denote the number
of homomorphisms  of $F$ into $G$ (adjacency-preserving maps from
$V(F)$ to $V(G)$), $\inj(F,G)$, the number of injective homomorphisms
of $F$ into $G$, and $\ind(F,G)$, the number of embedding of $F$ into
$G$ as an induced subgraph.

These quantities are closely related:
\[
\inj(F,G)=\sum_{F'\supseteq F}\ind(F',G),
\]
where $F'$ ranges over all graphs obtained from $F$ by adding edges,
and
\[
\hom(F,G)=\sum_{F''}\inj(F'',G),
\]
where $F''$ ranges over all graphs obtained from $F$ by identifying
nodes. Conversely, $\ind$ can be expressed by $\inj$, which in turn
can be expressed by $\hom$ using inclusion-exclusion.

This definition can be extended to the case when $G$ has nodeweights
$\alpha_v$ and edgeweights $\beta_{uv}$:
\[
\hom(F,G)=\sum_{\varphi:~V(F)\to V(G)} \prod_{u\in V(F)}
\alpha_{\varphi(u)}(G) \prod_{uv\in E(F)}
\beta_{\varphi(u),\varphi(v)}(G).
\]

We often normalize these homomorphism numbers, and consider the {\it
homomorphism densities}
\[
t(F,G)=\frac{\hom(F,G)}{|V(G)|^{|V(F)|}},
\]
which is the probability that a random map of $V(F)$ into $V(G)$ is a
homomorphism. We can define similarly
\begin{equation}\label{EQ:TINJ}
t_\inj(F,G)=\frac{\inj(F,G)}{n(n-1)\cdots(n-k+1)}
\end{equation}
and
\begin{equation}\label{EQ:TIND}
t_\ind(F,G)=\frac{\ind(F,G)}{n(n-1)\cdots(n-k+1)}.
\end{equation}
We have
\begin{equation}\label{EQ:TINJ-TIND}
t_\inj(F,G)=\sum_{F'\supseteq F}t_\ind(F',G)
\end{equation}
and the inversion formula
\begin{equation}\label{EQ:TIND-SZITA}
t_\ind(F,G)=\sum_{F'\supseteq F}(-1)^{|E(F')\setminus
E(F)|}t_\inj(F',G).
\end{equation}
For $\hom$ and $\inj$ the relationship is not so simple due to the
different normalization, but recalling that we are interested in
large graphs $G$, the following fact is usually enough to go between
them:
\begin{equation}\label{EQ:INJ-T}
t_\inj(F,G)-t(F,G)= O(\frac1{|V(G)|}).
\end{equation}

We note that $t_\ind(F,G)$ is the probability that sampling $V(F)$
nodes of $G$, we see the graph $F$. So it follows that (for very
large graphs, up to the error in \eqref{EQ:INJ-T}) subgraph sampling
provides the same information as any of the homomorphism densities
$t,t_\inj,t_\ind$.

\subsubsection{Spectra}

Homomorphisms of ``small'' graphs into $G$ are related to sampling,
as mentioned earlier. There are less obvious applications of these
numbers.

\begin{example}\label{CYCLES}
If $C_k$ denote the cycle on $k$ nodes, then $\hom(C_k,G)$ is the
trace of the $k$-th power of the adjacency matrix of the graph $G$.
In other words,
\[
\hom(C_k,G)=\sum_{i=1}^n \lambda_i^k,
\]
where $\lambda_1,\dots,\lambda_n$ are the eigenvalues of the
adjacency matrix of $G$. From here, eigenvalues with large absolute
value can be recovered. For example, $\hom(C_{2k},G)^{1/(2k)}$ tends
to the largest eigenvalue of $G$ as $k\to\infty$.
\end{example}

\subsection{Homomorphisms to the right}

\subsubsection{Colorings and independent sets}

Several important graph parameters can be expressed in terms of
homomorphisms into fixed ``small'' graphs.

\begin{example}\label{COLOR}
If $K_q$ denotes the complete graph with $q$ nodes (no loops), then
$\hom(G,K_q)$ is the number of colorings of the graph $G$ with $q$
colors, satisfying the usual condition that adjacent nodes must get
different colors.
\end{example}

\begin{example}\label{INDEP}
Let $H$ be obtained from $K_2$ by adding a loop at one of the nodes.
Then $\hom(G,H)$ is the number of independent sets of nodes in $G$.
\end{example}

\subsubsection{Multicuts}\label{MULTIWAY}

An important graph parameter is the {\it maximum cut} $\Maxcut(G)$,
the maximum number of edges between a set $S\subseteq V(G)$ of nodes
and its complement. While finding minimum cuts is perhaps more
natural, the maximum cut problem comes up when we want to approximate
general graphs by bipartite graphs, in computing ground states in
statistical physics (see next section), and in many other
applications. For our purposes, it will be more convenient to
consider the {\it normalized maximum cut}, defined by
\[
\maxcut(G)=\frac{\Maxcut(G))}{|V|^2}= \max_{S\subseteq V}
\frac{e_G(S,V\setminus S)}{|V|^2}
\]
(here $e_G(X,Y)$ denotes the number of edges in $G$ connecting node
sets $X$ and $Y$).

The following easy fact relates maximum cuts and homomorphism
numbers. Let $H$ be the weighted graph on $\{1,2\}$ with nodeweights
and edgeweights $1$ except for the non-loop edge, which has weight
$2$. Then we have the trivial inequalities
\[
2^{\Maxcut(G})\le \hom(G,H) \le 2^{|V(G)|}2^{\Maxcut(G)},
\]
which upon taking the logarithm and dividing by $|V(G)|^2$ becomes
\begin{equation}\label{EQ:CUT-HOM}
\maxcut(G)\le \frac{\log_2 \hom(G,H)}{|V(G)|^2} \le \maxcut(G)
+\frac{1}{|V(G)|}.
\end{equation}
So the homomorphism number into this simple 2-node graph determines
$\maxcut(G)$ asymptotically.

An important extension of the maximum cut problem involves partitions
into $q\ge 1$ classes instead of $2$. Instead of just counting edges
between different classes, we specify in advance numbers $\beta_{ij}$
($i,j\in[q]$) such that $\beta_{ij}=\beta_{ji}$. We define the
maximum multicut density (with the target weights $\beta_{ij}$) as
\[
\mmcut(G,\beta)= \max \frac1{|V(G)|^2}\sum_{i,j} \beta_{ij}
e_G(S_i,S_j),
\]
where the maximum is taken over all partitions $\{S_1,\dots,S_q\}$ of
$V(G)$.

A further important extension is to fix the proportion into which the
cut separates the node set. For example, the ``maximum bisection
problem'' asks for the maximum size of a cut that separates the nodes
into two equal parts (we allow a difference of $1$ if the number of
nodes is even). More precisely, we can formulate the {\it restricted
multicut problem} as follows. We specify (in addition to the
$\beta_{ij}$) numbers $\alpha_1,\dots, \alpha_q>0$ with
$\alpha_1+\cdots+\alpha_q=1$. It is convenient to consider the
parameters $\alpha_i$ and $\beta_{ij}$ as the nodeweights and edge
weights of a weighted graph $H$ with $V(H)=[q]$. Then are interested
in
\begin{equation}\label{EQ:EEGH}
\EE(G,H)=\max \frac1{|V(G)|^2}\sum_{i,j} \beta_{ij} e_G(S_i,S_j),
\end{equation}
where $\{S_1,\dots,S_q\}$ ranges over all partitions of $V(G)$ such
that
\begin{equation}\label{EQ:REST-PART}
\bigl||S_i|-\alpha_i|V(G)|\bigr|< 1\qquad (i=1,\dots,q).
\end{equation}
(This can be defined for all graphs $H$ with positive nodeweights, by
scaling the nodeweights so that they sum to $1$.)

The following extension of \eqref{EQ:CUT-HOM} is easy to prove: for
$H$ fixed and $|V(G)|\to\infty$,
\begin{equation}\label{EQ:MCUT-HOM}
\frac{\log_2 \hom(G,H)}{|V(G)|^2} = \mmcut(G,\beta)
+O(\frac{1}{|V(G)|}).
\end{equation}
(Note that $\log_2 \hom(G,H)/|V(G)|^2$ is asymptotically independent
of the node weights of $H$.)

The restricted maximum multicut problem is also related to counting
homomorphisms, but the relationship is a little more complicated. Let
$G$ be a (very large) simple graph and $H$, a weighted graph with
V(H)=[q]. In the definition of $t(G,H)$ we considered random maps
$V(G)\to V(H)$, where the image of each node is chosen independently
from the distribution on $V(H)$ defined by the node weights. For most
of these random maps $\varphi$, $|\varphi^{-1}(i)|\approx
\alpha_i(H)|V(G)||$ by the law of large numbers. It turns out that
often it is advantageous to restrict ourselves to maps that are
"typical" in this sense. More precisely, let $S(G,H)$ denote the set
of those maps $\varphi:~V(G)\to V(H)$ for which
$\bigl||\varphi^{-1}(i)|-\alpha_i|V(G)|\bigr|< 1$ for all $i\in
V(H)$. Using this notation, we can write
\[
\rmcut(G,H)=\max_{\varphi\in S(G,H)}\sum_{u,v\in V(G)}
\beta_{\varphi(u),\varphi(v)}.
\]
Let $\widetilde{H}$ be the weighted graph in which the edge weights
are $\widetilde{\beta}_{ij}=\exp(\beta_{ij})$ instead of
$\beta_{ij}$. If we define
\[
\hom^*(G,\widetilde{H})=\sum_{\varphi\in S(G,H)} \prod_{uv\in E(G)}
\widetilde{\beta}_{\varphi(u),\varphi(v)},
\]
then the following inequality analogous to \eqref{EQ:MCUT-HOM} holds
for $|V(G)|\to\infty$:
\begin{equation}\label{EQ:RMCUT-HOM}
\rmcut(G,H) = \frac{\log\hom^*(G,\widetilde{H})}{|V(G)|^2}+
O(\frac1{|V(G)|}).
\end{equation}

\subsubsection{Statistical physics}\label{STATPHYS}

Graph homomorphism functions can be used to express partition
functions of various statistical physical models. Two basic types of
such models are ``hard-core'' and ``soft-core''.

To describe an example of a hard-core model, let $G$ be an $n\times
n$ grid, and suppose that every node of $G$ (every ``site'') can be
in one of two states, ``UP'' or ``DOWN''. The properties of the
system are such that no two adjacent sites can be ``UP''. A
``configuration'' is a valid assignment of states to each node. The
number of configurations is the number of independent sets of nodes
in $G$, which in turn can be expressed as the number of homomorphisms
of $G$ into the graph $H$ consisting of two nodes, "UP" and "DOWN",
connected by an edge, and with an additional loop at "DOWN".

In a soft-core spin model the sites are again nodes of a graph $G$,
which can be in one of $q$ possible states. For any two states $i$
and $j$, we specify an ``energy of interaction'' in the form of a
real number $J_{ij}$. A given configuration (assignment of states) is
given by a map $\varphi: V(G)\to [q]$, and its ``energy density'' is
expressed as
\begin{equation}\label{EE}
\EE_{\varphi}=\frac 2{|V(G)|^2} \sum_{uv\in E(G)}
J_{\varphi(u),\varphi(v)},
\end{equation}
From this, one defines the partition function as
\begin{equation}\label{ZGJ}
Z(G,J)=\sum_{\varphi: V(G)\to [q]} \exp(-\EE_\varphi).
\end{equation}
Another important quantity is the {\it ground state energy}
\begin{equation}\label{EGJ}
\EE(G,J)= \min_{\varphi: V(G)\to [q]}\EE_\varphi.
\end{equation}
Note that both of these quantities are familiar: if we take
$\beta=-J$, then $\EE(G,J)=-\rmcut(G,\beta)$, and if we take
$\beta_{ij}=\exp(J_{ij})$, then $Z(G,J)=\hom(G,\exp(\beta))$. Even
restricted multiway cuts correspond to a quantity studied in
statistical physics: it is called {\it microcanonical ground state
energy} there.

The above definitions don't work well for dense graphs $G$: as
remarked after \eqref{EQ:MCUT-HOM}, the numbers $\log_2
\hom(G,H)/|V(G)|^2$ are essentially independent of the node weights
of $H$, so we loose information here. In the mean-field theory, we
define the mean field partition function of a simple graph $G$ by
\begin{equation}
\label{ZDEF} Z(G,J)=\sum_{\varphi: V(G)\to [q]}
e^{-|V(G)|\EE_\varphi}.
\end{equation}
The {\it free energy} is defined by
\begin{equation}\label{EQ:FREE-EN}
\FF(G,H)=-\frac{\ln Z(G,H)}{|V(G)|}.
\end{equation}
Note that the normalization is different from \eqref{ZGJ} in the
exponent and therefore we only divide by $|V(G)|$ (as opposed to
\eqref{EQ:MCUT-HOM}).

For more about this connection, we refer to \cite{BCLSV2}.

\subsection{Homomorphisms densities in the sparse case}

The best analogue for sparse graphs of the homomorphism density
$t(F,G)$ is
\begin{equation}\label{EQ:SPARSE-S}
s(F,G)=\frac{\hom(F,G)}{|V(G)|},
\end{equation}
which we consider for connected graphs $F$. We can interpret this
number as follows. For $u\in V(F)$ and $v\in V(G)$, let $\hom_{v\to
u}(F,G)$ denote the number of homomorphisms $\varphi$ of $F$ into $G$
with $\varphi(u)=v$. Now we fix any node $u$ of $F$ and select a
uniform random node $v$ of $G$. Then $s(F,G)$ is the expectation of
$\hom_{v\to u}(F,G)$. We can interpret
\[
s_\inj(F,G)=\frac{\inj(F,G)}{|V(G)|},\qquad
s_\ind(F,G)=\frac{\ind(F,G)}{|V(G)|}
\]
similarly.

\begin{remark}\label{REM:MAGNI-HOM}
For bounded degree graphs the order of magnitude of $\hom(F,G)$
(where $F$ is fixed and $V(G)$ tends to infinity) is $|V(G)|^{c(F)}$,
where $c(F)$ is the number of connected components of $F$. But since
$\hom(F,G)$ is multiplicative over the connected components of $F$,
we don't loose any information if we restrict the definition $s(F,G)$
to connected graphs $F$.
\end{remark}

The sparse homomorphism densities \eqref{EQ:SPARSE-S} contain the
same information as the distribution of neighborhood samples. The
proof of this is a bit trickier here than in the dense case.

From the interpretation of $s(F,G)$ given above, we see that it can
be obtained as the expectation of the number of $\hom_{u\to
v}s(F,\mathbf{B})$, where $\mathbf{B}$ is a random ball from the
neighborhood sample distribution $\rho_{G,r}$, with center $v$ and
radius $r=|V(F)$.

To compute the neighborhood sample distributions from the quantities
$s(F,G)$, we first express the quantities $s_\inj(F,G)$ via
inclusion-exclusion. By a similar argument, we can express the
quantities $s_\ind(F,G)$.

Next, we consider graphs $F$ together with maps
$\delta:~V(F)\to\{0,\dots,d\}$, and we determine the numbers
\[
s_\ind(F,\delta,G)=\frac{\ind(F,\delta,G)}{|V(G)|},
\]
where $\ind(F,\delta,G)$ is the number injections $\varphi:~V(F)\to
V(G)$ which embed $F$ in $G$ as an induced subgraph, so that the
degree of $\varphi(v)$ is $\delta(v)$. This is again done by an
inclusion-exclusion argument.

Given a ball $B$ of radius $r$, the fraction of nodes $v\in V(G)$ for
which $B(v,r)=B$ is $\sum_\delta \ind(B,\delta,G)$, where the
summation extends over all functions $\delta$ which assigns to each
node of $B$ at distance $<r$ from the root its degree in $B$. This
proves that homomorphism densities and neighborhood sampling are
equivalent.

\subsection{Characterizing homomorphism numbers}

Multigraph parameters of the form $\hom(\cdot,H)$, where $H$ is a
weighted graph, were characterized by Freedman, Lov\'asz and
Schrijver \cite{FLS}.

\begin{theorem}\label{THM:FLS}
Let $f$ be a graph parameter defined on multigraphs without loops.
Then $f$ is equal to $\hom(.,H)$ for some weighted graph $H$ on $q$
nodes if and only if it is reflection positive and $\rk(M(f,k))\le
q^k$ for all $k$.
\end{theorem}

Several improvements and versions of this result have been obtained.
It is shown in \cite{LSz5} that it is enough to assume the rank
condition for $k\le 2$. Analogous characterizations can be given for
graph parameters of the form $\hom(\cdot,H)$ where the nodeweights in
$H$ are all $1$ \cite{Sch3}, and where $H$ is an unweighted graph
without multiple edges (but with loops allowed) \cite{LSch2}. There
is also an analogous (dual) characterization of graph parameters of
the form $\hom(F,.)$, defined on simple graph with loops, where $F$
is also a simple graph with loops \cite{LSch2}. These results can be
extended to directed graphs, hypergraphs, semigroups, and indeed, to
all categories satisfying reasonable conditions \cite{LSch3}.

The two conditions on connection matrices in the theorem have
interesting uses of their own.

\subsubsection{Reflection positivity and extremal graph theory}

Theorem \ref{THM:GRAPHON-CHAR} will give a number of equivalent
(cryptographic) descriptions of limit objects of growing graph
sequences, and it can be used to characterize all reflection positive
graph parameters, see Corollary \ref{COR:GRAPHON-CHAR}.

Reflection positivity implies a number of very useful relations
between the densities of various subgraphs in a given graph, which in
turn can be used to prove results in extremal graph theory. We will
illustrate this in Section \ref{EXT-GRAPHONS}.

We'll return to applications of reflection positivity of connection
matrices in the context of continuous generalizations of graphs
(Section \ref{EXT-GRAPHONS}) and in extremal graph theory (Section
\ref{EXT-GRAPHONS}).

\subsubsection{Finite connection rank}

The finiteness of the rank of connection matrices is also
interesting. One reason to be interested in this question is the fact
that such a graph parameter can be evaluated in polynomial time for
graphs with bounded treewidth \cite{L2}.

There are several examples of graph parameters with finite connection
rank \cite{L1}: the number of perfect matchings, the number of all
matchings, the number of Hamiltonian cycles, any evaluation of the
Tutte polynomial.

A challenging problem is to determine all graph parameters for which
all the connection matrices have finite rank. Homomorphism functions
$\hom(.,H)$ are examples for every weighted graph $H$ (here the
nodeweights and edgeweights can be negative). Dual homomorphism
densities $\hom(F,.)$ also have finite connection rank. Every
evaluation of the Tutte polynomial is a further example.

Very recently Godlin and Makowski proved that {\it all graph
parameters which are evaluations of graph polynomials definable in
Monadic Second Order Logic have finite connection rank.} This result
can be used mostly as a tool to prove that certain properties are not
definable this way.

Further variants of this problem ask for the characterization of
graph parameters with exponentially bounded connection rank, or
polynomially bounded connection rank.

\subsection{Graph algebras}\label{ALGEBRA}

A {\it quantum graph} is defined as a formal linear combination of a
finite number of graphs with real coefficients. For every quantum
graph $x$, let $N(x)$ be the maximum number of nodes in the graphs
occurring in $x$ with nonzero coefficient. The definition of
$\hom(F,G)$ and $t(F,G)$ extends to quantum graphs linearly: if
$f=\sum_{i=1}^n \lambda_i F_i$ and $g=\sum_{j=1}^m \mu_j G_j$, then
we define
\[
\hom(f,g)=\sum_{i=1}^n \sum_{j=1}^m \lambda_i\mu_j \hom(F_i,G_j).
\]
Quantum graphs are useful in expressing various combinatorial
situations. For example, for any graph $F$ we define
\begin{equation}\label{EQ:SZITA}
\widehat{F}=\sum_{F': V(F')=V(F)\atop E(F')\supseteq E(F)}
(-1)^{|E(F')|} F'.
\end{equation}
Then $t(\widehat{F},G)$ is just the probability that a random map
$V(F)\to V(G)$ preserves adjacency as well as non-adjacency.

Let $f$ be any graph parameter and fix an integer $k\ge 0$. Let
$\QQ_k$ denote the (infinite dimensional) vector space of all
$k$-labeled quantum graphs. We can turn $\QQ_k$ into an algebra by
using $F_1F_2$ introduced above as the product of two generators, and
then extending this multiplication to the other elements linearly.
Clearly $\QQ_k$ is associative and commutative. The graph $O_k$ on
$k$ nodes with no edges is the multiplicative unit in $\QQ_k$. If all
nodes of $F$ are labeled, then both $F$ and the quantum graph
$\widehat{F}$ introduced above (keeping the node labels) are
idempotent: $F^2=F$ and $\widehat{F}^2=\widehat{F}$.

Every graph parameter $f$ can be extended linearly to quantum graphs,
and defines an inner product on $\QQ_k$ by
\begin{equation}\label{EQ:INNER}
\langle x,y\rangle := f(xy).
\end{equation}
This means that our graph algebra is a {\it Frobenius algebra} (see
\cite{Kock}). This inner product has nice properties, for example
\begin{equation}\label{EQ:ASSOC}
\langle x,yz\rangle=\langle xy,z\rangle.
\end{equation}

Let $\NN_k(f)$ denote the kernel of this inner product, i.e.,
\[
\NN_k(f):=\{x\in\QQ_k:~f(xy)=0~\forall y\in\QQ_k\}.
\]
Then we can define the factor algebra
\[
\QQ_k/f:=\QQ_k/\NN_k(f).
\]

\begin{example}\label{EXA:MATCH}
As an example, consider the number $\pm(G)$ of perfect matchings in
the graph $G$. It is a basic property of this value that subdividing
an edge by two nodes does not change it. This can be expressed as
$P_4-P_2\in\NN_2(\pm)$, where $P_k$ denotes the paths with $k$ nodes,
of which the two endnodes are labeled.
\end{example}

We can introduce a third ``product'': the tensor product $G\otimes H$
of a $k$-labeled graph $G$ and an $l$-labeled graph $H$ is defined as
the $(k+l)$-labeled graph obtained as the disjoint union of $G$ and
$H$, where the labels in $H$ are increased by $k$. If $k=l=0$, then
the tensor product is the same as the product in the algebra $\QQ_k$.

The parameter $f$ is reflection positive if and only if the inner
product \eqref{EQ:INNER} is positive semidefinite on $\QQ_k$;
equivalently, positive definite on $\QQ_k/f$, so it turns $\QQ_k/f$
into a Hilbert space. In fact, the factor algebra $\QQ_k/f$ is a
finite dimensional commutative $*$-algebra, which has both a
commutative and associative product and a positive definite inner
product, related by $\langle x,yz\rangle=\langle xy,z\rangle$.

The dimension of $\QQ_k/f$ is the rank of the connection matrix. If
this rank is a finite number $m$ and the parameter is reflection
positive, it follows that $\QQ_k/f$ is isomorphic $\R^m$ endowed with
the coordinate-wise product and the usual inner product.

There are many algebraically interesting connections between these
algebras, for example, there is an embedding given by the tensor
product
\begin{equation}\label{EQ:SUPERADD}
\QQ_k/f~\otimes~ \QQ_l/f \hookrightarrow \QQ_{k+l}/f,
\end{equation}
which shows that $\dim(\QQ_k/f)$ is a superadditive function of $k$.

This nice algebraic structure can be exploited in various ways
\cite{FLS,L2,LSos,LSz2}. Let us sketch the proof of Theorem
\ref{THM:FLS} in an (easier) special case: when there is no
degeneracy in the sense that the embedding in \eqref{EQ:SUPERADD} is
an isomorphism (this is in fact the generic case, which occurs
whenever $f=\hom(.,H)$, where $H$ has no ``twin'' nodes nor any
nontrivial automorphism). So we have $\dim(\QQ_k/f)=q^k$ for all $k$.

Let $p_1,\dots,p_q$ be the basis of $\QQ_1/f$ consisting of
idempotents (corresponding to the standard basis vectors in $\R^q$).
Define $p_\varphi = p_{\varphi(1)}\otimes \cdots\otimes
p_{\varphi(k)}$ for all $\varphi:~[k]\to[q]$, then the $k$-labeled
quantum graphs $p_\varphi$ form a basis of $\QQ_k/f$ consisting of
idempotents.

We can define a weighted complete graph $H$ on $[q]$ as follows: let
$\alpha_i=f(p_i)$ and define $\beta_{ij}$ by expressing the graph
$k_2$ (a single edge with both nodes labeled) in the idempotent
basis:
\[
k_2=\sum_{i,j\in[q]} \beta_{ij}(p_i\otimes p_j)
\]
This defines nodeweights $\alpha_i$ and edgeweights $\beta_{ij}$ for
$H$. The nodeweights are positive, since
\[
\alpha_i=f(p_i)=f(p_i^2)>0.
\]
The definition of the $\beta_{ij}$ implies that
\begin{equation}\label{EQ:K2B}
k_2(p_i\otimes p_j)=\beta_{ij}(p_i\otimes p_j).
\end{equation}
We claim that the weighted graph $H$ obtained this way satisfies
$f(G)=\hom(G,H)$ for every multigraph $G$. Indeed, we may assume that
$V(G)=[k]$ and all nodes of $G$ are labeled. Then we can write
\[
G=\prod_{uv\in E(G)} K_{uv},
\]
where $K_{uv}$ consists of $k$ labeled nodes and a single edge
connecting $u$ and $v$. Equation \eqref{EQ:K2B} implies that
\[
p_\varphi K_{uv} = \beta_{\varphi(u)\varphi(v)} p_\varphi.
\]
Using \eqref{EQ:ASSOC} repeatedly, we get
\[
G=\Bigl(\sum_{\varphi:\,[k]\to[q]} p_\varphi\Bigr)G =
\sum_{\varphi:\,[k]\to[q]} p_\varphi \prod_{uv\in
E(G)}K_{uv}=\sum_{\varphi:\,[k]\to[q]} \prod_{uv\in E(G)}
\beta_{\varphi(u)\varphi(v)} p_\varphi,
\]
and so
\[
f(G) =\sum_{\varphi:\,[k]\to[q]} \prod_{uv\in E(G)}
\beta_{\varphi(u)\varphi(v)}\prod_{u\in V(G)} \alpha_{\varphi(u)} =
\hom(G,H).
\]

\section{Graph-like structures on probability spaces}

The aim of this section is to introduce certain analytic objects,
which will serve as limit objects for graph sequences, separately in
the dense and sparse case. It is an interesting feature of these
structures that they have come up in different studies.

In the dense case, several versions of these objects turn out to be
equivalent; graphons are very simple objects (2-variable measurable
functions), but they turn out to be equivalent, among others, to
exchangeable random variables.

In the bounded degree case, several related, but non-equivalent
notions have been proposed, at least one of which (graphings) is also
known from group theory.

\subsection{Graphons}\label{GRAPHON}

Let $\WW$ denote the space of all bounded symmetric measurable
functions $W:~[0,1]^2\to\R$ (i.e., $W(x,y)=W(y,x)$ for all
$x,y\in[0,1]$). Let $\WW_0$ denote the set of all functions $W\in\WW$
such that $0\le W\le 1$.

A function $W\in\WW$ is called a {\it stepfunction}, if there is a
partition $S_1\cup\dots\cup S_k$ of $[0,1]$ into measurable sets such
that $W$ is constant on every product set $S_i\times S_j$. The number
$k$ is the {\it number of steps} of $W$.

For every weighted graph $G$, we define a stepfunction $W_G\in\WW_0$
as follows. Let $V(G)=[n]$. Split $[0,1]$ into $n$ intervals
$J_1,\dots,J_n$ of length $\lambda(J_i)=\alpha_i/\alpha_G$. For $x\in
J_i$ and $y\in J_j$, let
\[
W_G(x,y)=\beta_{ij}(G).
\]

Let $W\in\WW$ and let $\varphi:~[0,1]\to[0,1]$ be a measure
preserving map. We can define another function $W^\varphi$ by
\[
W^\varphi(x,y)=W(\varphi(x),\varphi(y)).
\]
From the point of view of using these functions as continuous
analogues of graphs, the functions $W$ and $W^\varphi$ are not
essentially different (they are related like two isomorphic graphs in
which the nodes are labeled differently). One has to be a little
careful though, because measure preserving maps are not necessarily
invertible, and so the relationship between $W$ and $W^{\varphi}$ is
not symmetric. We call two graphons $W$ and $W'$ {\it weakly
isomorphic}, if there is a third graphon $U$ and measure preserving
maps $\varphi,\varphi':~[0,1]\to[0,1]$ such that $W=U^\varphi$ and
$W'=U^{\varphi'}$ almost everywhere. It is not hard to show that weak
isomorphism is an equivalence relation.

Equivalence classes of functions in $\WW_0$ under weak isomorphism
are called {\it graphons}. (Sometimes we call a function $W\in\WW_0$
a graphon; by analogy with graphs, these functions could be called
``labeled graphons''.)

\subsubsection{Homomorphisms into graphons and from graphons}

Counting homomorphism into graphs extends to counting homomorphism
into graphons in the following sense: For every $W\in\WW$ and simple
graph $F=(V,E)$, define
\[
t(F,W)=\int_{[0,1]^V} \prod_{ij\in E} W(x_i,x_j)\,\prod_{i\in V} dx_i
\]
Then it is easy to verify that for every graph $G$,
\begin{equation}\label{EQ:TWG}
t(F,G)=t(F,W_G).
\end{equation}
Of the two modified versions of homomorphism densities
\eqref{EQ:TINJ} and \eqref{EQ:TIND}, the former has not significance
in this context since a random assignment $i\mapsto x_i$ ($i\in V(F),
x_i\in[0,1]$ is injective with probability $1$. But the induced
subgraph density is worth defining, and in fact it can be expressed
as
\begin{equation}\label{EQ:TIND-W}
t_{\ind}(F,W)=\int_{[0,1]^V} \prod_{ij\in E} W(x_i,x_j)\prod_{ij\in
\binom{V}{2}\setminus E} (1-W(x_i,x_j))\,\prod_{i\in V} dx_i.
\end{equation}
We have then
\begin{equation}\label{EQ:TINDWG}
t_\ind(F,G)=t_\ind(F,W_G),
\end{equation}
and the inclusion-exclusion formula \eqref{EQ:TIND-SZITA} follows by
expanding the parentheses in the integrand \eqref{EQ:TIND-W}.

Borgs, Chayes and Lov\'asz \cite{BCL} proved that the homomorphism
densities determine the graphon:

\begin{theorem}\label{THM:LIMUNIQUE}
Two graphons are weakly isomorphic if and only if $t(F,W)=t(F,W')$
for every simple graph $F$.
\end{theorem}

A natural idea of the proof of this theorem would be to bring every
graphon to a ``canonical form'', so that weakly isomorphic graphons
would have identical canonical forms. In the case of functions in a
single variable, a canonical form that works in many situations can
be obtained through ``monotonization'': for every bounded real
function on $[0,1]$ there is an unique monotone increasing
left-continuous function on $[0,1]$ that has, among others, the same
moments. For graphons this does not seem to be doable, but the proof
of Theorem \ref{THM:LIMUNIQUE} goes by constructing, for every
graphon $W$, a ``canonical ensemble'': a probability distribution on
graphons on the same canonical $\sigma$-algebra and weakly isomorphic
to $W$, such that two graphons are isomorphic if and only if their
ensembles can be coupled so that corresponding graphons are
identical.

Alternate proofs of Theorem \ref{THM:LIMUNIQUE} have been given by
Diaconis and Janson \cite{DiJa} using the theory of exchangeable
random variables, and by Bollob\'as and Riordan \cite{BoRi} combining
Theorem \ref{LEFT-CLOSE} below with measure-theoretic arguments.

There is probably no good way to define homomorphism numbers from
graphons into graphs or into other graphons. The parameters related
to such homomorphisms that extend naturally to graphons are defined
by maximization, like the normalized maximum cut, and more generally,
restricted maximum multiway cuts. Let $H$ be a weighted graph with
$V(H)=[q]$ and $W$, a graphon. Then we can define
\[
\EE(W,H)=\sup_{S_i} \sum_{i,j\in V(H)} \beta_{ij} \int_{S_i\times
S_j} W(x,y)\,dx\,dy,
\]
where $\{S_1,\dots,S_q\}$ ranges over all partitions of $[0,1]$ into
measurable sets with $\lambda(S_i)=\alpha_i(H)$. This quantity does
not exactly extend $\EE(G,H)$ as defined in \eqref{EQ:EEGH}, but the
error is small: it was proved in \cite{BCLSV2} that for a fixed
weighted graph $H$,
\begin{equation}\label{EQ:EEGHWH}
\EE(G,H)-\EE(W_G,H)=O\Bigl(\frac1{|V(G)|}\Bigr)\qquad
(|V(G)|\to\infty).
\end{equation}

\subsubsection{$W$-random graphs}\label{WRAND}

A graphon $W$ gives rise to a way of generating random graphs that
are more general than the Erd\H{o}s--R\'enyi graphs. This
construction was introduced by Lov\'asz and Szegedy \cite{LSz1} and
Bollob\'as, Janson and Riordan \cite{BJR}.

Given a graphon $W$ and an integer $n>0$, we can generate a random
graph $\Ge(n,W)$ on node set $[n]$ as follows: We generate $n$
independent numbers $X_1,\dots,X_n$ from the uniform distribution on
$[0,1]$, and then connect nodes $i$ and $j$ with probability
$W(X_i,X_j)$, making an independent decision for distinct pairs
$(i,j)$.

As a special case, if $W$ is the identically $p$ function, we get
``ordinary'' random graphs $\Ge(n,p)$.

We can extend this construction to generating a countable random
graph $\Ge(W)$ on $\N$: We generate an infinite sequence
$X_1,X_2,\dots$ of uniformly distributed random points from $[0,1]$,
and (as before) connect nodes $i$ and $j$ with probability
$W(X_i,X_j)$.

Graphons will come up in several ways in our discussions. In Theorem
\ref{THM:GRAPHON-CHAR} we will collect the many disguises in which
they occur.

\subsection{Graphings}\label{GRAPHING}

\subsubsection{Measure preserving graphs}

Let $G$ be a graph with node set $[0,1]$, with all degrees bounded by
$d$. We call $G$ {\it measurable}, if for every (Lebesgue) measurable
set $B$ the neighborhood $N(B)$ in $G$ is also measurable.

For every set $A\subseteq[0,1]$ and $x\in[0,1]$, let $d_A(x)$ denote
the number of neighbors of $x$ in $B$. One can prove using the
measurability of $G$ that $d_A(x)$ is a measurable function of $x$.
We say that $G$ is {\it measure preserving}, if it is measurable and
for any two measurable sets $A,B$,
\begin{equation}\label{EQ:UNIMOD}
\int_A d_B(x)\,dx =\int_B d_A(x)\,dx.
\end{equation}

Assuming that this relation holds, we can define a measure $\mu$ on
the Borel sets of $[0,1]^2$ by $\mu(A\times B)=\int_A d_B(x)\,dx$.
This measure is concentrated on the set of edges (which can be
considered as a subset of $[0,1]^2$). Furthermore, the marginals of
$\mu$ are absolutely continuous with respect to the Lebesgue measure,
and their Radon-Nikodym derivative is the degree function.

In every measure preserving graph $G$, we can define the density
$s(F,G)$ of a graph $F$. Indeed, let us recall that $s(F,G)$ is the
expectation of $\hom_{v\to u}(F,G)$, where $v$ is a fixed node of $F$
and $u$ is a random node of $G$. Since we have a probability
distribution on $V(G)$, and $\hom_{v\to u}(F,G)$ is a bounded
measurable function of $u$, this definition carries over verbatim.

Similarly, we can talk about the neighborhood distributions
$\rho_{G,m}$ in a measure preserving graph.

\subsubsection{Graphings}

Let $A_1,\dots,A_d,B_1,\dots,B_d$ be measurable subsets of $[0,1]$,
and let $\varphi_i:~A_i\to B_i$ be bijective measure preserving maps.
The tuple $H=([0,1],\varphi_1,\dots,\varphi_d)$ is called a {\it
graphing} (see \cite{Gab,KM}). From every graphing $H$ we get a
directed graph $\overrightarrow{G}$ on $[0,1]$ by connecting $x$ and
$y$ in $[0,1]$ if there is an $i$ such that $y=\varphi_i(x)$. The
edges of this digraph are colored with $d$ colors in such a way that
each color-class defines a measure preserving bijection between two
subsets of $[0,1]$.

Forgetting the orientation and the edge-colors of this digraph, we
get a measure preserving graph with degrees bounded by $2d$. A
measure preserving graph with its edges colored and oriented so that
each color defines a measure preserving bijection is equivalent to a
graphing.

It would be perhaps more natural to assume that the maps
$\varphi_1,\dots,\varphi_d$ are involutions, in which case we get an
undirected graph, and we can extend the $\varphi_i$ to measure
preserving involutions $[0,1]\to[0,1]$. It is true that for every
graphing there is such an involutive graphing defining the same
measure preserving graph; but the number of maps may become much
larger.

Every measure preserving graph arises from a graphing:

\begin{theorem}\label{THM:GRAPH2GRAPHING}
Let $G$ be a measure preserving graph with degrees bounded by $d$.
Then there is a graphing $H=([0,1],\varphi_1,\dots,\varphi_r)$, where
$r\le d^2$, such that the underlying graph is $G$.
\end{theorem}

One way of looking at a representation of a measure preserving graph
as a graphing is that it provides a {\it certificate} that the graph
is measure preserving. The graphing representing a given measure
preserving graph may not be unique.

Theorem \ref{THM:GRAPH2GRAPHING} can be viewed as a measure
preserving graph version of Shannon's Theorem, which asserts that the
edges of a multigraph with maximum degree $d$ can be colored by
$3d/2$ colors. (For simple graphs, Vizing's Theorem gives the better
bound of $d+1$.) The bound $d^2$ is probably not optimal in the
measure preserving version either.

We will talk about $s(F,H)$ if $F$ is a (finite) graph and $H$ is a
graphing. This will mean simply $s(F,G)$, where $G$ is the underlying
measure preserving graph.

We note that both in measure preserving graphs and graphings, we
could replace the probability space $[0,1]$ by any other standard
probability space, but this would not lead to any gain in generality.
However, in some cases the presentation of the measure preserving
graph or graphing is more natural on other probability spaces.

\subsubsection{Random countable rooted graphs}

Measure preserving graphs are also related to certain probability
distributions on rooted countable graphs, introduced by Benjamini and
Schramm \cite{BSch}.

Let $G$ be a measure preserving graph and choose a uniform random
point $x\in[0,1]$. The connected component $G_x$ of $G$ containing
$x$ is a countable graph with degrees bounded by $d$, and with a
``root'' node $x$.

Let $\Gb_d$ denote the set of connected countable graphs with all
degrees bounded by $d$, rooted at a node. Let $\AA_d$ denote the
$\sigma$-algebra on $\Gb_d$ generated by subsets obtained by fixing a
finite neighborhood of the root. The map $x\mapsto G_x$ is measurable
as a map $[0,1]\to(G_d,\AA_d)$, and thus every measure preserving
graph $G$ defines a probability distribution $\pi$ on
$(\Gb_d,\AA_d)$.

Condition \ref{EQ:UNIMOD} implies the following property of the
measure $\pi$. Selecting a rooted graph $G$ from $\pi$ and then
selecting a uniform random edge from the root, we get a probability
distribution $\pi^*$ on the set $\Gb'_d$ of rooted graphs in $\Gb_d$
with an edge (the ``root edge'') from the root also specified. We say
that $\pi$ is {\it unimodular}, if the map $\Gb'_d\to\Gb'_d$ obtained
by shifting the root node to the other endnode of the root edge is
measure preserving with respect to $\pi$.

The measure on $\Gb_d$ obtained from a measure preserving graph is
unimodular. Vice versa, every such measure is obtained from a
graphing (and hence from a measure preserving graph; Elek
\cite{Elek1}).

\section{The cut-distance of two graphs}\label{DIST}

The definition of the distance of two arbitrary graphs is quite
involved, and we will approach the problem in steps: starting with
two graphs on the same node set, then moving to graphs with the same
number of nodes (but unrelated), then moving to the general case.

In this section we consider dense graphs. The definitions are of
course valid for all graphs, but they give a distance of $o(1)$
between two graphs with edge-density $o(1)$.

\subsection{Two graphs on the same set of nodes}

Let $G$ and $G'$ be two graphs with a common node set $[n]$. The
distance notion discussed here was initiated by Frieze and Kannan
\cite{FK}, and elaborated, e.g., in \cite{BCLSV1}. For an unweighted
graph $G=(V,E)$ and sets $S,T\subseteq V$, let $e_G(S,T)$ denote the
number of edges in $G$ with one endnode in $S$ and the other in $T$
(the endnodes may also belong to $S\cap T$; so $e_G(S,S)$ is twice
the number of edges spanned by $S$). For two graphs $G$ and $G'$ on
the same node set $[n]$, we define their {\it cut distance} by
\[
d_\square(G,G')=\frac{1}{n^2}\max_{S,T\subseteq V(G)}
|e_G(S,T)-e_{G'}(S,T)|.
\]
Note that we are dividing by $n^2$ and not by $|S|\times |T|$, which
would look more natural. However, dividing by $|S|\times |T|$ would
emphasize small sets too much, and the maximum would be attained when
$|S|=|T|-1$. With our definition, the contribution of a pair $S,T$ is
at most $|T|\cdot|S|/n^2$ (for simple graphs).

It is easy to see that $d_\square(G,G')\le d_1(G,G')$, and in general
the two sides are quite different. For example, if $\Ge$ and $\Ge'$
are two independent random graphs on $[n]$ with edge probability
$1/2$, then with large probability $d_\square(G,G')=O(1/\sqrt{n})$.

\subsection{Two graphs with the same number of nodes}

If $G$ and $G'$ are unlabeled unweighted graphs on different node
sets but of the same cardinality $n$, then we define their distance
by
\begin{equation}\label{EQ:DIST-UNW-UNL}
\hat\delta_\square(G,G')=\min_{\tilde{G},\tilde{G}'}
d_\square(\tilde{G},\tilde{G}'),
\end{equation}
where $\tilde{G}$ and $\tilde{G}'$ range over all labelings of $G$
and $G'$ by $1,\dots,n$, respectively. (The hat above the $\delta$
indicates that the ``ultimate'' definition will be somewhat
different.)

\subsection{Two arbitrary graphs}

Let $G=(V,E)$ and $G'=(V',E')$ be two graphs with (say) $V=[n]$ and
$V'=[n']$. To define their distance, we need a graph operation: for
every graph $G$ and positive integer $m$, let $G(m)$ denote the graph
obtained from $G$ by replacing each node of $G$ by $m$ nodes, where
two new nodes are connected if and only if their predecessors were.

We can use the distance $\hat\delta_\square$ to define the distance
\[
\delta_\square(G,G')=\lim_{k\to\infty}
\hat\delta_\square(G[kn'],G'[kn]).
\]
(Here $G(kn')$ and $G'(kn)$ have the same number of nodes.)

A more complicated but ``finite'' definition of the same quantity can
be given as follows. A {\it fractional overlay} of $G$ and $G'$ is a
nonnegative $n\times n'$ matrix $X$ such that $\sum_{u=1}^{n'}
X_{iu}=\frac{1}{n}$ and $\sum_{i=1}^{n} X_{iu}=\frac{1}{n'}$. If
$n=n'$ and $\sigma:~V\to V'$ is a bijection, then $X_{iu}=\frac1n
\one_{\sigma(i)=u}$ is a fractional overlay (which in this case is an
honest-to-good overlay). We denote by $\XX(G,G')$ the set of all
fractional overlays.

For a matrix $M$, let $\Sigma(M)$ denote the sum of its entries. Then
the distance of the two graphs can be described by the following
optimization problem:
\begin{equation}\label{DIST-W-UNL}
\delta_\square(G,G')=\min_{X\in \XX(G,G')} \max_{Y,Z\subseteq V\times V'}
\Bigl|\sum_{iu\in Y,~jv\in Z\atop ij\in E}X_{iu}X_{jv}-\sum_{iu\in
Y,~jv\in Z\atop uv\in E'}X_{iu}X_{jv}\Bigr|.
\end{equation}

To illuminate this definition a little, we can think of a fractional
overlay as a coupling of the uniform distribution on $V(G)$ with the
uniform distribution on $V(G')$: it gives a probability distribution
$\chi$ on $V(G)\times V(G')$ whose marginals are uniform. Select two
pairs $(i,u)$ and $(j,v)$ from the distribution $\chi$. Then the
first sum in \eqref{DIST-W-UNL} is the probability that ``$iu\in Y$
and $jv\in Z$ and $ij\in E$'', and the second sum is the probability
that ``$iu\in Y$ and $jv\in Z$ and $uv$ is an edge''. Thus
\eqref{DIST-W-UNL} expresses some form of correlation between $ij$
being an edge and $uv$ being an edge.

One word of warning: $\delta_\square$ is only a pseudometric, not a
true metric, because $\delta_\square(G,G')$ may be zero for different
graphs $G$ and $G'$. This is the case e.g. if $G'=G(k)$ for some $k$.

Definition \eqref{DIST-W-UNL} can be extended to weighted graphs, but
instead of going through the hairy formulas, we postpone this to the
next section.

We conclude with a problem for which only partial results are
available. If $G$ and $G'$ have the same number of nodes, then the
definition of $\delta_\square$ does not give back
$\hat\delta_\square$. It was proved in \cite{BCLSV1} that
\begin{equation}\label{EQ:SQ-HATSQ}
\delta_\square(G,G')\le\hat\delta_\square(G,G')\le
32\delta_\square(G,G)^{1/67}.
\end{equation}
This is a rather weak result, its significance being that
$\delta_\square$ and $\hat\delta_\square$ define the same Cauchy
sequences. Alon (unpublished) proved that
\begin{equation}\label{EQ:SQ-HATSQ2}
\hat\delta_\square(G,G')\le (1+o(1))\delta_\square(G,G)
\end{equation}
if $|V(G)|=|V(G')|\to\infty$. We conjecture:

\begin{conj}\label{DDHAT}
For any two graphs $G$ and $G'$ on $n$ nodes,
$\hat\delta_\square(G,G') \le 2\delta_\square(G,G')$.
\end{conj}

An analogous result for the edit distance was proved by Pikhurko
\cite{Pik}.

\subsection{Distance of graphons}

This notion of distance extends to graphons as follows (and it is
perhaps more natural in that context). We consider on $\WW$ the {\it
cut norm}
\[
\|W\|_\square = \sup_{S,T\subseteq[0,1]}\Bigl|\int_{S\times T}
W(x,y)\,dx\,dy\Bigr|
\]
where the supremum is taken over all measurable subsets $S$ and $T$.
It is sometimes convenient the use the corresponding metric
$d_1(U,W)=\|U-W\}_\square$. We define the {\it cut distance}
\[
\delta_\square(U,W)=\inf_{\varphi} d_\square(U,W^\varphi),
\]
where $\varphi$ ranges over all invertible measure preserving maps
from $[0,1]\to[0,1]$, and $W^\varphi(x,y)=W(\varphi(x),\varphi(y))$.

The distance $\delta_\square$ of graphons is only a pseudometric,
since different graphons can have distance zero. This happens
precisely when they are weakly isomorphic.

If $G$ and $G'$ are weighted graphs, then we have
\begin{equation}\label{EQ:DEF-DELTA-SQ}
\delta_\square(G,G')=\delta_\square(W_G,W_{G'}).
\end{equation}
This could serve as a more natural (but not combinatorial) definition
of the distance of two graphs, and we will use it to define the
distance of two weighted graphs. Let $K$ denote the graph with a
single node of weight $1$, endowed with a loop with weight $1/2$.
Then for a random graph $\Ge=\Ge(n,1/2)$, we have
$\delta_\square(\Ge,K)=O(1/\sqrt{n})$ with large probability.

Going into all the complications with using the cut norm and then
minimizing over measure preserving transformations is justified by
the following important fact.

\begin{theorem}\label{THM:COMPACT}
The pseudometric space $(\WW_0,\delta_\square)$ is compact.
\end{theorem}

The proof depends on Szemer\'edi partitions, to be discussed in
section \ref{SZEMEREDI}.

Convergence in the $\|.\|_\square$ norm is stronger than
weak-$*$-convergence. To be more precise, if $\|W_n-W\|_\square\to 0$
($n\to\infty$), then it follows immediately from the definition that
\begin{equation}\label{EQ:WEAK-ST}
\int_{S\times T} W_n\to \int_{S\times T} W,
\end{equation}
and hence by standard arguments we get that
\begin{equation}\label{EQ:WEAK-ST2}
\int_{[0,1]^2}U\cdot W_n\to \int_{[0,1]^2} U\cdot W
\end{equation}
for every integrable function $U$. However, weak-$*$-convergence is
not equivalent of convergence in the $\|.\|_\square$ norm; a
counterexample can be obtained e.g. from Example \ref{EXA:PFX-ATT}
(see \cite{BCLSV3}).

Similar construction can be applied to other norms, e.g., from the
$L_1$-norm
\[
\|W\|_1=\int_{[0,1]^2} |W(x,y)|\,dx
\]
we get
\[
d_1(U,W)=\|U-W^\varphi\|_1\qquad \text{and}\qquad
\delta_1(U,W)=\inf_\varphi \delta_1(U,W).
\]

\section{Szemer\'edi partitions}\label{SZEMEREDI}

One of the most important tools in understanding large dense graphs
is the Regularity Lemma of Szemer\'edi \cite{Szem1,Szem2} and its
extensions. This lemma has many interesting connections to other
areas of mathematics, including analysis \cite{LSz3,BoNi} and
information theory \cite{Tao1}. It also has weaker (but more
effective) and stronger versions. Here we survey as much as we need
from this rich theory.

\subsection{$\eps$-regular bipartite graphs and the original lemma}

For a graph $G=(V,E)$ and for $X,Y\subseteq V$, let $e_G(X,Y)$ denote
the number of edges with one endnode in $X$ and another in $Y$; edges
with both endnodes in $X\cap Y$ are counted twice. We denote by
$d_G(X,Y)=\frac{e_G(X,Y)}{|X|\cdot|Y|}$ the density of edges between
$X$ and $Y$. If $X$ and $Y$ are disjoint, we denote by $G[X,Y]$ the
bipartite graph on $X\cup Y$ obtained by keeping just those edges of
$G$ that connect $X$ and $Y$.

Let $\PP=\{V_1,\dots,V_k\}$ be a partition of $V$. We say that $\PP$
is an {\it equipartition} if $\lfloor|V|/k\rfloor \le |V_i| \le
\lceil |V_i|/k\rceil$ for all $1\le i\le k$. We define the weighted
graph $G_\PP$ on $V$ by taking the complete graph and weighting its
edge $uv$ by $d_G(V_i,V_j)$ if $u\in V_i$ and $v\in V_j$.

The Regularity Lemma says, roughly speaking, that every graph has a
partition $\PP$ into a ``small'' number of classes such that $G_\PP$
is ``close'' to $G$. There are (non-equivalent) forms of this lemma,
depending on how we measure the error.

Let $G$ be a bipartite graph $G$ with bipartition $\{U,W\}$. On the
average, we expect that for $X\subseteq U$ and $Y\subseteq W$,
\[
e_G(X,Y)\approx d_G(X,Y) |X|\cdot|Y|.
\]
For two arbitrary subsets of the nodes, $e_G(X,Y)$ may be very far
from this ``expected value'', but if $G$ is a random graph, then,
however, it will be close; random graphs are very ``homogeneous'' in
this respect. We say that $G$ is {\it $\eps$-regular}, if
\begin{equation}\label{EQ:E-REGULAR}
\left|\frac{e_G(X,Y)}{|X|\cdot|Y|}-d\right| \le\eps
\end{equation}
holds for all subsets $X\subseteq U$ and $Y\subseteq W$ such that
$|X|>\eps|U|$ and $|Y|>\eps |W|$.

Notice that we could not require condition \eqref{EQ:E-REGULAR} to
hold for small $X$ and $Y$: for example, if both have one element,
then the quotient $e_G(X,Y)/(|X|\cdot|Y|)$ is either 0 or 1. However,
we could replace it by the condition
\begin{equation}\label{EQ:E-REGULAR2}
\bigl|e_G(X,Y)-d|X|\cdot|Y|\bigr| \le\eps|U|\cdot|W|
\end{equation}
for all $Y\subseteq U$ and $Y\subseteq W$. Indeed,
\eqref{EQ:E-REGULAR} implies \eqref{EQ:E-REGULAR2} for $|X|>\eps|U|$
and $|Y|>\eps |W|$, while if e.g. $|X|\le\eps|U|$, then
$e_G(X,Y)\le\eps|U|\cdot|W|$ and $d|X|\cdot|Y|| \le\eps|U|\cdot|W|$,
so \eqref{EQ:E-REGULAR2} holds trivially. Conversely, if
\eqref{EQ:E-REGULAR2} holds with $\eps$ replaced by $\eps^3$, then
\[
\left|\frac{e_G(X,Y)}{|X|\cdot|Y|}-d\right|\le\frac{\eps^3
|U|\cdot|W|}{|X|\cdot|Y|}<\eps
\]
if $|X|>\eps|U|$ and $|Y|>\eps |W|$.

With these definitions, the Regularity Lemma can be stated as
follows:

\begin{lemma}[Szemer\'edi Regularity Lemma, usual form]\label{SZEM-1}
For every $\eps>0$ there is a $k=k(\eps)$ such that every graph
$G=(V,E)$ on at least $k$ nodes has an equipartition
$\{V_1,\dots,V_k\}$ $(1/\eps\le k\le k(\eps))$ such that for all but
$\eps k^2$ pairs of indices $1\le i<j\le k$, the bipartite graph
$G[V_i,V_j]$ is $\eps$-regular.
\end{lemma}

One feature of the Regularity Lemma, which unfortunately forbids
practical applications, is that $k(\eps)$ is very large: the best
proof gives a tower of height about $1/\eps^2$, and unfortunately
this is not far from the truth, as was shown by Gowers \cite{Gow1}.

\subsection{Weak Regularity Lemma and distance of graphs}

A version with a weaker conclusion but with a more reasonable error
bound was proved by Frieze and Kannan \cite{FK}.

\begin{lemma}[Weak Regularity Lemma]\label{W-SZEM}
For every $k\ge 1$ and every graph $G=(V,E)$, $V$ has a partition
$\PP$ into $k$ classes such that
\[
d_\square(G,G_\PP)\le \frac{2}{\sqrt{\log k}}.
\]
\end{lemma}

Note that we do not require here that $\PP$ be an equipartition; it
is not hard to see that this version implies that there is also an
equipartition with similar property, just we have to increase the
error bound to $4/\sqrt{\log k}$.

To see the connection with the original lemma, we note that if $G$ is
an {\it $\eps$-regular} bipartite graph say in the sense of
\eqref{EQ:E-REGULAR2}, and $H$ is the weighted complete bipartite
graph with the same bipartition $\{U,W\}$ and with edge weights $d$,
then \eqref{EQ:E-REGULAR2} says that $d_\square(G,H)\le\eps$. Hence
if $\PP$ is a Szemer\'edi partition in the sense of Lemma
\ref{SZEM-1}, then the distance between the bipartite subgraph of $G$
induced by $V_i$ and $V_j$, and the corresponding weighted bipartite
subgraph of $G_\PP$, is at most $\eps$ for all but $\eps k^2$ pairs
$(i,j)$, and at most $1$ for the remaining $\eps k^2$ pairs. It is
easy to see that this implies that the distance between $G$ and
$G_\PP$ is at most $\eps$. So the partition in Lemma \ref{W-SZEM} has
indeed weaker properties than the partition in Lemma \ref{SZEM-1}. Of
course, this is compensated for by the relatively decent number of
partition classes.

If we keep in $G_\PP$ an edge with weight $p$ with probability $p$
and delete it with probability $1-p$, then we get a random graph $H$,
and it is easy to see that with large probability
$d_\square(G_\PP,H)\le \frac{10}{\sqrt{|V(G)|}}$. This implies the
following version of the Weak Regularity Lemma:

\begin{lemma}\label{LEM:W-SZEM-SIMPLE}
For every $k\ge 1$ and graph $G$, there is a graph $H$ with
$k$ nodes such that
\[
\delta_\square(G,H)\le \frac{10}{\sqrt{\log k}}.
\]
\end{lemma}

\subsection{Strong Regularity Lemma and compactness}

Other versions of the Regularity Lemma strengthen, rather than
weaken, the conclusion (of course, at the cost of replacing the tower
function by an even more formidable value). Such a ``super-strong''
Regularity Lemma was proved by Alon, Fisher, Krivelevich and Szegedy
\cite{AFKS}. We state the following equivalent version from
\cite{LSz3}.

\begin{lemma}[Strong Regularity Lemma]\label{ST-SZEM}
For every sequence $(\eps_0,\eps_1,...)$ of positive numbers there is
a positive integer $k_0$ such that for every graph $G=(V,E)$, there
is a graph $G'$ on $V$, and $V$ has a partition $\PP$ into $k\le k_0$
classes such that
\begin{equation}\label{EQ:SRL}
d_1(G,G')\le\eps_0 \qquad\text{and}\qquad d_\square(G',G'_\PP)\le
\eps_k.
\end{equation}
\end{lemma}

\noindent Note that the first inequality involves the normalized edit
distance, and so it is stronger than a similar condition with the cut
distance would be. The second error bound $\eps_k$ in \eqref{EQ:SRL}
can be thought of very small. If we choose $\eps_k=\eps_0$, we get
the Frieze--Kannan version \ref{W-SZEM} (with $\eps=2\eps_0$).
Choosing $\eps_k=\eps_0/k^2$, the partition obtained satisfies the
requirements of the original Regularity Lemma \ref{SZEM-1}.

The strong version itself follows rather easily from the compactness
of the space $(\WW_0,\delta_\square)$ (Theorem \ref{THM:COMPACT});
see \cite{LSz1} for details.

\subsection{Partitions into sets with small diameter}

\subsubsection{Small diameter sets and regularity}

We can equip every graph $G=(V,E)$ with a metric as follows. Let $A$
be the adjacency matrix of $G$. We define the {\it similarity
distance} of two nodes $i,j\in V$ as the $\ell_1$ distance of the
corresponding rows of $A^2$ (squaring the matrix seems unnatural, but
it is crucial; it turns out to get rid of random fluctuations). The
following was proved (in somewhat different form) in \cite{LSz3}.

\begin{theorem}\label{LEM:METRIC-SZEM}
Let $G$ be a graph and let $\PP=\{V_1,\dots,V_k\}$ be a partition of
$V$.

\smallskip

{\rm (a)} If $d_\square(G,G_\PP)=\eps$, then there is a set
$S\subseteq V$ with $|S|\le 8\sqrt{\eps}|V|$ such that for each
partition class, $V_i\setminus S$ has diameter at most $8\sqrt{\eps}$
in the $d_2$ metric.

\smallskip

{\rm (b)} If there is a set $S\subseteq V$ with $|S|\le
\delta|V|$ such that for each partition class,
$V_i\setminus S$ has diameter at most $\delta$ in the $d_2$ metric,
then $d_\square(G,G_\PP)\le 24\delta$.
\end{theorem}

Theorem \ref{LEM:METRIC-SZEM} suggests to define the {\it dimension}
of a family $\GG$ of graphs as the infimum of real numbers $d>0$ for
which the following holds: for every $\eps>0$ and $G\in\GG$ the node
set of $G$ can be partitioned into a set of at most $\eps|V(G)|$
nodes and into at most $\eps^{-d}$ sets of diameter at most $\eps$.
(This number can be infinite.) In the cases when the graphs have a
natural dimensionality, this dimension tends to give the right value.
For example, let $G$ be obtained by selecting $n$ random points on
the $d$-dimensional unit sphere, and connecting two of these points
$x$ and $y$ with a probability $W(x,y)$, which is a continuous
function of $x$ and $y$. With probability $1$, this sequence has
dimension $\Theta(d)$.

\subsubsection{Computational applications}\label{SZEM-ALG}

As an easy application of Theorem \ref{LEM:METRIC-SZEM}, we give an
algorithm to compute a weak Szemer\'edi partition in a huge graph.
Our goal is to illustrate how an algorithm works in the pure sampling
model, as well as in what form the result can be returned. This way
of presenting the output of an algorithm for a large graph was
proposed by Frieze and Kannan \cite{FK}.

We start with an auxiliary algorithm that computes (approximately)
the $d_2$ distance of two nodes.

\begin{alg}\label{ALG:D2-AUX}\strut

{\bf Input:} A graph $G$ given by an sampling oracle, two nodes
$u,v\in V$, and an error bound $\eps>0$.

\smallskip

{\bf Output:} A number $D_2(u,v)\ge 0$ such that with probability at
least $1-\eps$,
\[
D_2(u,v)-\eps\le d_2(u,v) \le D_2(u,v)+\eps.
\]
\end{alg}

To see how this can be done, we rewrite the definition of the $d_2$
distance as follows. For $x,y\in V(G)$, let $a(x,y)$ be the
corresponding entry of the adjacency matrix of $G$: this is $1$ if
they are adjacent and $0$ otherwise. Define
\[
a_2(x,y) = \E_z a(x,z)a(y,z),
\]
where $z$ is a uniform random node in $V$; this is the corresponding
entry of the square of the adjacency matrix, normalized by $|V(G)|$.
Finally, let
\[
d_2(x,y) = \E_z(|a_2(x,z)-a_2(y,z)|),
\]
where again $z$ is a uniform random node in $V$. Drawing a
sufficiently large sample (depending on $\eps$), these expectations
can be approximated by averaging.

Algorithm \ref{ALG:D2-AUX} enables us to encode a partition of $V(G)$
as a subset $R\subseteq V(G)$: for each $r\in R$, we define the
partition class $V_r$ as the set of nodes $u\in V$ such that the node
in $R$ closest to $u$ is $r$. Ties will be broken arbitrarily, and
nodes to which there are several ``almost closest'' nodes may be
misclassified, but this is the best one can hope for. To formalize,

\begin{alg}\label{ALG:SZEM-AUX}\strut

{\bf Input:} A graph $G$ given by an sampling oracle, a subset
$R\subseteq V(G)$, a node $u\in V$, and an error bound $\eps>0$.

\smallskip

{\bf Output:} An $r\in R$ such that with probability at least
$1-\eps$, $d_2(u,r)\le (1+\eps)d_2(u,r)$.
\end{alg}

The way this second algorithm works is that it uses Algorithm
\ref{ALG:D2-AUX} to compute (approximately) the distances $d_2(u,r)$,
$r\in R$, and returns the node $r\in R$ that it finds closest to $u$.
Borrowing a phrase from geometry, we compute the Voronoi cells of the
set $R$.

Using this encoding of the partition, the following algorithm
computes a weak Szemer\'edi partition.

\begin{alg}\label{ALG:SZEM}\strut

{\bf Input:} A graph $G$ given by an sampling oracle, and an error
bound $\eps$.

\smallskip

{\bf Output:} A set $R\subseteq V(G)$ with $|R|\le 2^{2/\eps^2}$ such
that, with probability at least $1-\eps$, $d_2(u,R)\le \eps$ for all
but an $\eps$ fraction of the nodes $u$.
\end{alg}

The set $R$ is grown step by step, starting with the empty set. At
each step, a new uniform random node $w$ of $G$ is generated, and the
approximate distances $D_2(u,v)$ are computed for all $r\in R$ with
error less than $\eps/|R|$. If all of these are larger than $\eps/2$,
$w$ is added to $R$. Else, $w$ is thrown out and a new random node is
generated. If $R$ is not increased in $1/\eps^2$ steps, the algorithm
halts.

It is clear that if more than an $\eps$ fraction of the nodes are
farther than $\eps$ from $R$, then in $1/\eps^2$ iterations we are
very likely to sample one of these, and then with large probability
we get the distances right and so we increase $R$.

Theorem \ref{LEM:METRIC-SZEM} says in this context that the partition
determined by Algorithms \ref{ALG:D2-AUX}--\ref{ALG:SZEM} satisfies
$d_\square(G,G_\PP)\le (4\eps)^{1/4}$ with large probability.

We conclude with an answer to Question 4 in Section \ref{QUEST}. For
the partition $\PP$ implicitly determined above, we can also compute
the edge densities between the partition classes, which we use to
weight the edges of the complete graph on $R$, so that we get a
weighted graph $H$. We find the maximum cut in $H$ by brute force, to
get a partition $R=R_1\cup R_2$. This gives an implicit definition of
a cut in $G$, where a node $u$ if put on the left side of the cut iff
$D_2(u,R_1)<D_2(u,R_2)$ for the approximate distances computed by
Algorithm \ref{ALG:D2-AUX}.

\subsection{Regularity Lemmas for bounded degree graphs?}

The Regularity Lemma as discussed above does not say anything for
non-dense graphs. Several extensions for this case are known
\cite{Koh,GS}, but they are meaningless for graphs that are very
sparse, in particular if they have bounded degree.

Is there a Regularity Lemma for graphs with bounded degree? There are
great difficulties here, but three results justify cautious optimism.

An observation of Alon (unpublished) implies that a weak analogue of
the Regularity Lemma, version \ref{LEM:W-SZEM-SIMPLE}, holds. Using
the sampling distance introduced in Section \ref{LOC-SAMP}, we can
state this as follows:

\begin{prop}\label{PROP:ALON-B-REG}
For every $d\ge 1$ and $\eps>0$ there is a $n=n(d,\eps)$ such that
for every graph $G$ with degrees bounded by $d$ there is a graph $H$
with degrees bounded by $d$ and $|V(H)|\le n$, such that $d_{\rm
sample}(G,H)\le \eps$.
\end{prop}

Unfortunately, no effective bound on $n$ follows from the proof. It
would be very interesting to give any explicit bound on the function
$n(d,\eps)$, or to give an algorithm to construct $H$ from $G$.
Ideally, one would like to design an algorithm that would work in the
sampling framework, similarly as the algorithm in Section
\ref{SZEM-ALG} works in the dense case.

It was proved recently by Elek and Lippner \cite{ElekLip}, and
independently by Angel and Szegedy \cite{ASz} that every graph with
degrees bounded by $d$ can be decomposed by deleting $\eps n$ edges
into ``highly homogeneous'' parts, where the number of these parts is
bounded by a function of $d$ and $\eps$. Unfortunately, the highly
homogeneous parts can still have a complex structure, but this may be
a first important step in the direction of finding an analogue of the
Regularity Lemma.

A third idea of decomposition is related to F{\o}lner sequences in
the theory of amenable groups, and is called {\it hyperfiniteness}
for general graph sequences \cite{Elek5,Schramm1}. A family $\GG$ of
graphs with bounded degree is called {\it hyperfinite}, if for every
$\eps>0$ there is a $k_\eps\ge 1$ such that from every graph
$G\in\GG$ we can delete $\eps|V(G)|$ edges so that every connected
component of the remaining graph has at most $k_\eps$ nodes. Schramm
\cite{Schramm1} showed that for a convergent graph sequence,
hyperfiniteness is reflected by the limit object.

A special case of a hyperfinite family is a family $\GG$ of graphs
with {\it subexponential growth}, familiar from group theory. This
property is defined by requiring that there is a function
$f:~\N\to\N$ such that $(\ln f(m))/m\to 0$ ($m\to\infty$), and for
any graph $G\in\GG$, any $v\in V(G)$ and any $m\in\N$, the number of
nodes in the $m$-neighborhood of $v$ is at most $f(m)$.

It is likely that large real-life networks can be thought of as
hyperfinite; on the other hand, hyperfinite families of graphs seem
to be much better behaved, and some of the theory of dense graph
sequences can be extended at least to this case.

\section{Convergence and limits I: the dense case}\label{CONV-LIM-DENSE}

\subsection{Subgraph sampling}\label{CONV-DENSE}

Recall that we can define a notion of convergence if we fix a {\it
sampling method}. For dense graphs, we use subgraph sampling: We
select uniformly a random $k$-element subset of $V(G)$, and return
the subgraph $G[k]$ induced by them. The probability that we see a
given graph $F$ is the quantity $t_\ind(F,G)$ introduced in
\eqref{EQ:TIND}. A sequence of graphs $(G_n)$ with
$|V(G_n)|\to\infty$ is {\it convergent} if the induced subgraph
densities $t_\ind(F,G_n)$ converge for all finite graphs $F$.

We use this sampling method for dense graphs (otherwise all these
densities tend to $0$).

Instead of the induced subgraph densities $t_\ind(F,G_n)$, we could
use the subgraph densities $t_\inj(F,G_n)$ or the homomorphism
densities $t(F,G_n)$. Indeed, the subgraph densities can be expressed
as linear combinations of induced subgraph densities and vice versa,
while the difference $t(F,G)-t_\inj(F,G)=O(1/|V(G)|)$, and so it
tends to $0$ if $|V(G)|\to\infty$.

We can extend this sampling procedure to graphons, and we get to the
construction of $W$-random graphs.

\subsection{Convergence in distance}

The definition of convergence can be reformulated using the notion of
sampling distance \ref{EQ:SAMPDIST}: a sequence $(G_n)$ of simple
graphs with $|V(G_n)|\to\infty$ is convergent if for every graph $F$,
$(t_\ind(F,G_n):~i=1,2,\dots)$ is a Cauchy sequence (equivalently,
$(t(F,G_n):~i=1,2,\dots)$ is a Cauchy sequence). This is equivalent
to saying that the graph sequence is Cauchy in the $d_{\rm sample}$
metric. The following theorem, which is one of the main results in
this theory, justifies the use of the cut metric $\delta_\square$.

\begin{theorem}\label{THM:CAUCHY}
A sequence $(G_n)$ of simple graphs ($|V(G_n)|\to\infty$) is
convergent if and only if it is a Cauchy sequence in the metric
$\delta_\square$.
\end{theorem}

A quantitative form of this equivalence is given by the following
theorem. Part (a) is a generalization of what is called the
``Counting Lemma'' in the theory of Szemer\'edi partitions; part (b)
may be called the ``Anti-counting'' lemma.

\begin{theorem}\label{LEFT-CLOSE}
Let $U,W\in\WW_0$.

\smallskip

{\rm (a)} For every simple finite graph $F$,
\[
|t(F,U)-t(F,W)|\le |E(F)|\cdot\delta_\square(U,W).
\]

{\rm (b)} Let $k$ be a positive integer, and assume that for every
simple graph $F$ on $k$ nodes, we have
\[
|t(F,U)-t(F,W)|\leq 2^{-{k^2}}.
\]
Then
\[
\delta_\square(U, W) \le \frac{20}{\sqrt{\log k}}.
\]
\end{theorem}

The proof of part (a) is quite simple; part (b) depends on the
sampling lemmas to be discussed in Section \ref{SEC:SAPLEMMAS}.

Theorem \ref{THM:CAUCHY} can be generalized to characterize
convergence in the space $\WW$:

\begin{theorem}\label{LEFTCONV=DELTACONV}
Let $(W_n)$ be a sequence of graphons in $\WW_0$ and let $W\in
\WW_0$. Then $t(F,W_n)$ converges for all finite simple graphs $F$ if
and only if $W_n$ is a Cauchy sequence in the $\delta_\square$
metric. Furthermore $t(F,W_n)\to t(F,W)$ for all finite simple graphs
$F$ if and only if $\delta_\square(W_n,W)\to 0$.
\end{theorem}

\subsection{Convergence from the right}

Convergence of a graph sequence can also be characterized in terms of
mappings ``to the right''. Several characterizations along these
lines were given in \cite{BCLSV2}; here we state one:

\begin{theorem}\label{THM:RIGHT-CONV}
Let $(G_n)$ be a sequence of simple graphs such that
$|V(G_n)|\to\infty$ as $n\to\infty$.  Then the sequence $(G_n)$ is
left-convergent if and only if the sequence $\EE(G_n,H)$ is
convergent for every weighted graph $H$.
\end{theorem}

\subsection{Sampling and distance}\label{SEC:SAPLEMMAS}

The proof of the results in the previous section depends on a couple
of probabilistic lemmas, which relate sampling to cut distance. The
first of these lemmas is due to Alon, Fernandez de la Vega, Kannan
and Karpinski \cite{AFKK}, with an improvement from \cite{BCLSV1}.
Its proof is quite involved. Its main implication is that the
$d_\square$-distance of two graphs $G$ and $H$ on the same set of
nodes can be estimated by sampling. It should be noted that the bound
given is quite sharp.

\begin{lemma}\label{LEM:SAMPLE1}
Let $k$ be a positive integer and let $G$ and $H$ be graphs with
$V(G)=V(H)$, $|V(G)|\ge k$ and edge weights in $[0,1]$. Let $S$ be
chosen uniformly from all subsets of $V(G)$ of size $k$. Then with
probability at least $1-2e^{-\sqrt k/8}$.
\[
\Bigl| d_\square(G[S],H[S])-d_\square(G,H)\Bigr| \le
\frac{10}{k^{1/4}}.
\]
\end{lemma}

The second lemma about sampling \cite{BCLSV1} shows that a sample is
close to the original graph with large probability. Note that here we
have to sue the $\delta_\square$ distance (since no overlaying is
given a priori), and also that the bound on the distance is much
weaker than in the previous lemma.

\begin{lemma}\label{LEM:SAMPLE2}
Let $k\ge 1$, and let $G$ be a simple graph on at least $k$ nodes. If
$S$ is a random subset of $V(G)$ of size $k$, then with probability
at least $1-2^{-k}$,
\[
\delta_\square(G, G[S]) \le \frac{10}{\sqrt{\log k}}.
\]
\end{lemma}

This lemma follows from Lemma \ref{LEM:SAMPLE1} and the Weak
Regularity Lemma \ref{W-SZEM}. Let us sketch this proof.

\begin{proof}
Fix some $m\ge 1$. By Lemma \ref{W-SZEM}, there is an equipartition
$\PP=\{V_1,\dots,V_m\}$ of $V(G)$ into $m$ classes such that
\[
d_\square(G,G_\PP)\le \frac{4}{\sqrt{\log m}}.
\]
Now let $S$ be a random $k$-subset. By Lemma \ref{LEM:SAMPLE1}, we
have
\[
|d_\square(G[S],G_\PP[S])-d_\square(G,G_\PP)|\le \frac{10}{k^{1/4}}
\]
with large probability. If $k$ is sufficiently large relative to $m$,
then every class $V_i$ will contain about $k/m$ nodes from $S$.
Indeed, a simple application of Chebyshev's Inequality gives that
with probability at least $3/4$,
\[
\Bigr||V_i\cap S|-\frac{k}{m}\Bigr| \le 2\sqrt{k}
\]
holds for all $i$.

Now blow up each node of $G_\PP[S]$ into $n/k$ twins to get a
weighted graph $G'$ (in notation: $G'=G_\PP[S](\frac nk)$). Then each
set $V_i\cap S$ is blown up into a set $V_i'$ of size
$\frac{k}{n}|V_i\cap S|\approx|V_i|=\frac{n}{m}$. In fact,
\[
\sum_{i=1}^m \bigl| |V_i'|-|V_i|\bigr| \le
2\sqrt{k}m\frac{n}{k}=\frac{2nm}{\sqrt{k}}.
\]
It follows that we can overlay $G'$ and $G_\PP$ so that corresponding
edges have the same weight except for edges inside the classes $V_i$
and edges incident with at most $\sum_{i=1}^m \bigl| |V_i'|-
|V_i|\bigr|$ nodes. This is only a fraction of $\frac{1}{m} +
\frac{4m}{\sqrt{k}}$ of all edges, which shows that
\[
\delta_\square(G_\PP,G_\PP[S])=\delta_\square(G_\PP,G')\le \frac{1}{m} +
\frac{4m}{\sqrt{k}}.
\]
Hence
\begin{align*}
\delta_\square(G,G[S])&\le \delta_\square(G,G_\PP)+\delta_\square(G_\PP,G_\PP[S])
+\delta_\square(G_\PP[S],G[S])\\
&\le \frac{4}{\sqrt{\log m}}+\frac{10}{k^{1/4}}+\Bigl( \frac{1}{m} +
\frac{4m}{\sqrt{k}}\Bigr).
\end{align*}
Choosing $m=k^{1/4}$, we get
\[
\delta_\square(G,G[S])\le \frac{8}{\sqrt{\log
k}}+\frac{10}{k^{1/4}}+\Bigl( \frac{1}{k^{1/4}} +
\frac{4}{k^{1/4}}\Bigr) < \frac{10}{\sqrt{\log k}}
\]
if $k$ is large enough.
\end{proof}

Both lemmas \ref{LEM:SAMPLE1} and \ref{LEM:SAMPLE2} extend to
graphons. We only formulate the second one, which can be stated in
terms of the $W$-random graphs $\Ge(k,W)$.

\begin{lemma}\label{LEM:SAMPLE2a}
Let $k\ge 1$, and let $W$ be a graphon. Then with probability at
least $1-2^{-k}$,
\[
\delta_\square(\Ge(k,W),W) \le \frac{11}{\sqrt{\log k}}.
\]
\end{lemma}

To illustrate how these lemmas fit in the proofs, let us first sketch
how Lemma \ref{LEM:SAMPLE2a} implies the ``anti-counting lemma''
(Theorem \ref{LEFT-CLOSE}(b)). Assume that $U,W\in\WW_0$ satisfy
\[
|t(F,U)-t(F,W)|\leq 2^{-{k^2}}
\]
for every graph $F$ with $k$ nodes. In terms of the $W$-random graphs
$\Ge(k,U)$ and $\Ge(k,W)$, this implies (by inclusion-exclusion) that
\[
\bigl| \Pr(\Ge(k,U)\cong F)-\Pr(\Ge(k,W)\cong F)\bigr| \le
2^{\binom{k}{2}}2^{-{k^2}},
\]
and hence
\[
\sum_F \bigl| \Pr(\Ge(k,U)\cong F)-\Pr(\Ge(k,W)\cong F)\bigr| \le
2^{k(k-1)}2^{-{k^2}}= 2^{-k}.
\]
This means that we can couple $\Ge(k,U)$ and $\Ge(k,W)$ so that
$G(k,U)\cong G(k,W)$ with probability at least $1-2^{-k}$. Lemma
\ref{LEM:SAMPLE2a} implies that with probability at least $1-2^{-k}$,
we have
\[
\delta_\square(U,\Ge(k,U))\le \frac{10}{\sqrt{\log k}},
\]
and similar assertion holds for $W$. Whenever all three happen, we
get
\[
\delta_\square(U,W)
\le\delta_\square(U,\Ge(k,U))+\delta_\square(\Ge(k,U),\Ge(k,W))
+\delta_\square(W,\Ge(k,W)) \le \frac{20}{\sqrt{\log k}}.
\]

\subsection{Dense limit}\label{LIM-DENSE}

The main motivation behind considering graphons is the following
theorem \cite{LSz1}:

\begin{theorem}\label{thm:LIMIT-EXIST}
For any convergent sequence $(G_n)$ of simple graphs there exists a
graphon $W$ such that $t(F,G_n)\to t(F,W)$ for every simple graph
$F$.
\end{theorem}

We say that this graphon $W$ is the {\it limit} of the graph
sequence, and write $G_n\to W$.

One might wonder if we really need complicated objects like
integrable functions to describe these limits; would perhaps
piecewise linear, or monotone, or continuous functions suffice? The
following two results tell us that (up to weak isomorphism) all
measurable functions are needed: every graphon $W$ can be obtained as
the limit of a sequence of simple graphs \cite{LSz1}, and the limit
is essentially unique \cite{BCL}.

\begin{theorem}\label{WRAND-CONV}
For any $W\in\WW_0$, the graph sequence $\Ge(n,W)$ converges to the
graphon $W$ with probability $1$.
\end{theorem}

On the other hand, Theorem \ref{THM:LIMUNIQUE} implies:

\begin{theorem}[\cite{BCL}]\label{thm:LIMIT-UNIQ}
The limit graphon of a convergent graph sequence is uniquely
determined up to weak isomorphism.
\end{theorem}

There are two quite different proofs of the (main) theorem
\ref{thm:LIMIT-EXIST}. The original one in \cite{LSz1} uses
Szemer\'edi partitions and the Martingale Convergence Theorem; a more
recent proof by Elek and Szegedy \cite{ESz} first constructs a
different limit object in the form of an uncountable graph by taking
the ultraproduct, and them obtains the graphon as an appropriate
projection of this (in terms of non-standard analysis, the graphon is
a non-standard Szemer\'edi partition of this graph on a non-standard
$[0,1]$ interval).

The first proof has the obvious advantage of being a constructive;
but the second proof is very general, it extends to hypergraphs and
many other structures, and leads to new understanding of the
Regularity Lemma for hypergraphs \cite{Gow1,Gow2,RS1} and its
consequences \cite{Tao2}.

Convergence to the limit object can also be characterized by the
distance function introduced above \cite{BCLSV1}:

\begin{theorem}\label{THM:CONV-METRIC}
For a sequence $(G_n)$ of graphs with $|V(G_n)|\to\infty$ and graphon
$W$, we have $G_n\to W$ if and only if $\delta_\square(W_{G_n},W)\to
0$.
\end{theorem}

Note that the function $W_{G_n}$ depends on the labeling of the nodes
of $G_n$ (the distance $\delta_\square(W_{G_n},W)$ does not, since
relabeling $G_n$ results in weak isomorphism of $W_{G_n}$). Choosing
the labeling appropriately, we can say more:

\begin{theorem}\label{THM:CONV-NORM}
For a sequence $(G_n)$ of graphs with $|V(G_n)|\to\infty$ and graphon
$W$, we have $G_n\to W$ if and only if the graphs $G_n$ can be
labeled so that $\|W_{G_n}-W\|_\square\to 0$.
\end{theorem}

\subsubsection{Equivalent descriptions of the limit}

A {\it random graph model} is a probability distribution on simple
graphs on $[n]$, for every $n\ge 1$, which is invariant under the
reordering of the nodes. In other words, it is a sequence of random
variables $\Ge_n$, whose values are simple graphs on $[n]$, and
isomorphic graphs have the same probability. We say that a random
graph model is {\it consistent} if deleting node $n$ from $\Ge_n$,
the distribution of the resulting graph is the same as the
distribution of $\Ge_{n-1}$. We say that the model is {\it local}, if
for every $1<k<n$, the subgraphs of $\Ge_n$ induced by $[k]$ and
$\{k+1,\dots,n\}$ are independent as random variables.

It is easy to see that for every graphon $W\in\WW_0$, $\Ge(n,W)$ is a
consistent and local random graph model.

A related notion is the following. Let $\GG$ be the set of graphs on
$\N$; we can think of $\GG$ as the product space $\{0,1\}^E$, where
$E=\binom{\N}{2}$ is the set of all (unordered) pairs of elements of
$\N$. This also equips $\GG$ with a $\sigma$-algebra. Let $\Sigma$ be
the group of permutations of $\N$, and let $\Sigma_2$ be the action
of $\Sigma$ on $E$. Recall that a probability measure $\pi$ on $\GG$
is called {\it ergodic} with respect to $\Sigma_2$ if it is invariant
under $\Sigma_2$ and $\GG$ has no measurable subset $\GG'$ with
$0<\pi(\GG')<1$ invariant under $\Sigma_2$.

It is easy to see for that every $W\in\WW$, the random graph $\Ge(W)$
defines a probability measure on $\GG$ invariant under $\Sigma_2$.
B.~Szegedy \cite{Sze2} showed that this measure is also ergodic.

After this preparation, we can formulate the theorem describing the
many notions equivalent to graphons.

\begin{theorem}\label{THM:GRAPHON-CHAR}
The following structures are cryptomorphic:

\smallskip

{\rm (a)} a graphon $W\in\WW_0$, up to weak isomorphism;

\smallskip

{\rm (b)} A graph parameter $f$ that is the limit of graph parameters
$t(.,G_n)$ for some convergent graph sequence $(G_n)$.

\smallskip

{\rm (c)} A multiplicative, reflection positive graph parameter $f$
satisfying $f(K_1)=1$,

\smallskip

{\rm (d)} a consistent local random graph model;

\smallskip

{\rm (e)} an ergodic measure on $\GG$ invariant under $\Sigma_2$.
\end{theorem}

The equivalences of these structures are mostly contained in results
mentioned previously. Let us sketch these constructions.

(a)$\to$(b): Every graphon $W\in\WW_0$ gives rise to the graph
parameter $t(.,W)$; furthermore, $W$ is the limit of a convergent
graph sequence $(G_n)$ (for example, of the sequence of $W$-random
graphs), and for this sequence $t(F,G_n)\to t(F,W)$ for all $F$.

(b)$\to$(c): If a graph parameter is the limit of graph parameters
$t(.,G_n)$, which satisfy the conditions in (c), then clearly so does
their limit.

(c)$\to$(d): In the special case when $f=t(.,G)$ is the probability
that a random map from $F$ to some graph $G$ is a homomorphism, we
can express the probability that a sample of $n$ points gives a given
graph $F_0$, by inclusion-exclusion in terms of the numbers $f(F)$.
We can apply the same formula to any graph parameter $f$ satisfying
(c), and get a probability distribution on $n$-point graphs (here the
conditions in (c) are used), which is a consistent local random graph
model.

(d)$\to$(a): Generating a random graph $\Ge_n$ from the consistent
local random graph model, it can be shown that we get a convergent
graph sequence with probability $1$, which tends to a graphon $W$.
For this graphon, $\Ge(n,W)$ gives back the random graph model we
started with.

(d)$\leftrightarrow$(e): It is easy to see that a consistent random
graph model is equivalent to a probability distribution on $\GG$
invariant under $\Sigma_2$. The proof that locality is equivalent to
ergodicity \cite{Sze2} is trickier and not given here.

\begin{corollary}\label{COR:GRAPHON-CHAR}
A graph parameter $f$ is reflection positive if and only if it is
either identically $0$, or there is a probability distribution $\rho$
on the Borel sets of $(\WW_0,\delta_\square)$ such that if
$\mathbf{W}$ denotes a random function from this distribution, then
\[
f(F)=\E t(F,\mathbf{W}).
\]
\end{corollary}

\subsubsection{Examples}\label{LIM-EXAMPLES}

We start with two easy examples.

\begin{example}\label{EXA:COMP-BIP}
{\it Complete bipartite graphs.} It is natural to guess, and easy to
prove, that complete bipartite graphs $K_{n,n}$ converge to the
function defined by $W(x,y)=1$ if $0\le x\le 1/2\le y\le 1$ or $0\le
y\le 1/2\le x\le 1$, and $W(x,y)=0$ otherwise.
\end{example}

\begin{example}\label{EXA:THRESH}
{\it Threshold graphs.} These graphs are defined on the set
$\{1,\dots,n\}$ by connecting $i$ and $j$ if and only if $i+j\le n$.
These graphs converge to the function defined by $W(x,y)=\one_{x+y\le
1}$.
\end{example}

\begin{example}\label{EXA:QUASIRAND}
A sequence of graphs tending to the identically-$p$ function is
exactly what we called a quasirandom sequence with density $p$.
\end{example}

Two examples of randomly growing graph sequences:

\begin{example}\label{EXA:UNIF-ATT}
{\it Randomly grown uniform attachment graph.} We start with a single
node. At the $n$-th iteration, a new node is born, and then every
pair of nonadjacent nodes is connected with probability $1/n$. We
call this graph sequence a {\it randomly grown uniform attachment
graph sequence}.

Let us do some simple calculations. After $n$ steps, let
$\{0,1,\dots,n-1\}$ be the nodes (born in this order). The
probability that nodes $i < j$ are not connected is
$\frac{j}{j+1}\cdot \frac{j+1}{j+2}\cdots\frac{n-1}{n}= \frac{j}{n}$.
These events are independent for all pairs $(i,j)$. From here, one
can easily figure out that the expected number of edges is
$(n^2-1)/6$.

To describe the limit function, note that the probability that nodes
$i$ and $j$ are not connected is $\max(i,j)/n$. If $i=xn$ and $j=yn$,
then this is $\max(x,y)$. Using that these events are independent, we
can prove that {\it the graph sequence $G^{\rm ua}_n$ tends to the
limit function $1-\max(x,y)$ with probability $1$.}
\end{example}

\begin{example}\label{EXA:PFX-ATT}
{\it Randomly grown prefix attachment graph.} In this construction,
it will be more convenient to label the nodes starting with $1$. At
the $n$-th iteration, a new node $n$ is born, a node $z\le n$ is
selected at random, and the new node is connected to nodes
$1,\dots,z-1$. We denote the $n$-th graph in the sequence by
$G_n^{\rm pfx}$, and call this graph sequence a {\it randomly grown
prefix attachment graph sequence}.

The expected number of edges is $n(n-1)/4$, and one can compute
subgraph densities with some effort to see that the sequence is
convergent with probability $1$. It is more difficult to figure out
the limit graphon.

We can try to proceed similarly as in the case of uniform attachment
graphs. The probability that $i$ and $j$ are connected is
$|j-i|/\max(i,j)$; if $i=xn$ and $j=yn$, then this is
$|x-y|/\max(x,y)$. Does this mean that the function
$U(x,y)=|x-y|/\max(x,y)$ is the limit? Surprisingly, the answer is
negative, which we can see by computing triangle densities.

The key to describe the limit is the remark at the end of Section
\ref{INTRO:LIMITS}, namely that instead of the uniform distribution
over the interval $[0,1]$, we can use other probability spaces. Let
us label a node born in step $k$, connected to $\{1,\dots,m\}$, by
the pair $(k/n, m/k)\in [0,1]\times[0,1]$. Then we can observe that
nodes with label $(x_1, y_1)$ and $(x_2, y_2)$ are connected if and
only if either $x_1< x_2y_2$ or $x_2< x_1y_1$.

From this observation one can prove that {\it the prefix attachment
graphs $G_n^{\rm pfx}$ converge, with probability 1, to the function
$W:~[0,1]^2\times [0,1]^2\to[0,1]$, given by}
\[
W_{\rm pfx}((x_1, y_1),(x_2, y_2))=
  \begin{cases}
    1, & \text{if $x_1< x_2y_2$ or $x_2< x_1y_1$}, \\
    0, & \text{otherwise}.
  \end{cases}
\]

This gives a nice and simple representation of the limit object with
the underlying probability space $[0,1]^2$ (with the uniform
measure). If we want a representation on $[0,1]$, we can map $[0,1]$
into $[0,1]^2$ by any measure preserving map $\varphi$; then $W_{\rm
pfx}^\varphi(x,y)=W^{\rm pfx}(\varphi(x),\varphi(y))$ gives a weakly
isomorphic graphon. This function is $0-1$ valued, but its support is
fractal-like.

It is interesting to note that the graphs $\Ge(n,W)$ form a different
growing sequence of random graphs tending to the same limit $W$ with
probability 1.
\end{example}

\subsection{Convergence from the right}\label{DENSE-CONV-R}

Paper \cite{BCLSV2} contains several conditions that characterize
convergent dense graph sequences in terms of homomorphisms ``to the
right'' (we have seen that these correspond to parameters with
meaning in statistical physics). We only state one of these, in our
terms:

\begin{theorem}\label{THM:RIGHTCONV}
Let $(G_n)$ be a sequence of graphs such that $|V(G_n)|\to\infty$ as
$n\to\infty$. Then the sequence $(G_n)$ is convergent if and only if
the restricted multicut densities $\rmcut(G_n,H)$ are convergent for
every weighted graph $H$.
\end{theorem}

By \eqref{EQ:RMCUT-HOM}, the value $\rmcut$ in this theorem could be
replaced by $\hom^*(G,H)$, and by our discussion in Section
\ref{STATPHYS}, we could talk about microcanonical ground state
energies instead of restricted multicuts.

\subsection{Limits of other dense combinatorial structures}

Limit objects can be defined for multigraphs, directed graphs,
colored graphs, hypergraphs etc. In many cases, like directed graphs
without parallel edges, or graphs with nodes colored with a fixed
number of colors, this can be done along the same lines as for simple
graphs.

But in other cases there are some surprises. For example, limits of
multigraphs with edge-multiplicities are not real valued functions,
but 2-variable functions whose values are random variables with
nonnegative integral values \cite{LSz5}. If $W$ is such a function,
we can generate a $W$-random multigraph by selecting $n$ independent
random points $X_1,\dots,X_n$ from the uniform distribution on
$[0,1]$, and then connecting nodes $i$ and $j$ with $W(X_i,X_j)$
parallel edges (which is a random integer).

The case of hypergraphs is much more interesting and important.
Formulating regularity lemmas and constructing limits of sequences of
$r$-uniform hypergraphs, where $r$ is fixed, is a highly nontrivial
task, but it is essentially solved now, thanks to the work of R\"odl
and Skokan and Gowers \cite{RS1,Gow}; see also \cite{Tao1,ESz}.

However, it seems that no good extension of the distance
$\delta_\square$ has been found to hypergraphs (just as for the
regularity lemma, the first natural guesses are wrong). Another open
question is to extend these results to nonuniform hypergraphs, with
unbounded edge-size.

The semidefiniteness conditions for homomorphism functions can be
extended to hypergraphs (see e.g. \cite{LSch1}). One area of
applications of these conditions is extremal graph theory, and it is
natural to ask if the semidefiniteness conditions can be useful in
extremal hypergraph theory, especially since extremal problems for
hypergraphs tend to be much harder than for graphs, and even basic
questions are unsolved.

\section{Convergence and limits II: bounded degree graphs}
\label{CONV-LIM-SPARSE}

\subsection{Neighborhood sampling}\label{CONV-SPARSE}

Recall the sampling process for bounded degree graphs: For a fixed
nonnegative integer $r$, we select uniformly a random node $v\in
V(G)$, and return the ball $B_G(v,r)$ with center $v$ and radius $r$
(i.e., the subgraph induced by those nodes that can be reached from
$v$ on a path of length $r$ or less). For a given rooted graph $F$,
we denote by $\rho_{G,r}(F)$ the probability that this sampling
method returns $F$ (with the root as the center). So $\rho_{G,r}(F)$
defines a probability distribution on rooted graphs $F$ with radius
at most $k$, which we denote by $\rho_{G,r}$.

We use this method if the degrees of nodes in $G$ are bounded by a
fixed number $d$; then the number of possible neighborhoods is
finite.

A sequence of graphs $(G_n)$ with degrees uniformly bounded by $d$
and $|V(G_n)|\to\infty$ is {\it convergent} (or more precisely
locally convergent) if the neighborhood densities $\rho_{G_n,r}(F)$
converge for all $r$ and all finite rooted graphs $F$.

Similarly as for the subgraph sampling, there are equivalent density
type parameters whose convergence could be used instead of the
neighborhood densities, for example, we could stipulate the
convergence of $s(F,G_n)$ for every connected graph $F$.

\subsection{Local (weak) limit}\label{LIM-SPARSE}

\subsubsection{Different forms}

A weakly convergent bounded degree graph sequence has several, not
quite equivalent limit objects, which we have introduced in Section
\ref{GRAPHING}. Part (a) of the following theorem is due to Benjamini
and Schramm \cite{BSch}; part (b) was formulated by R.~Kleinberg
(unpublished); part (c), which implies (b), is due to Elek
\cite{Elek1}.

\begin{theorem}\label{THM:WLIM}
Let $(G_n)$ be a locally convergent sequence of graphs with degrees
bounded by $d$. Then

\smallskip

{\rm (a)} There is a unique unimodular distribution $\tau$ on
countable rooted graphs with degrees bounded by $d$ such that
$\rho_{G_n,r}\to\rho_\tau$.

\smallskip

{\rm (b)} There is a measure preserving graph $G$ such that
$\rho_{G_n,r}\to \rho_{G,r}$ for every $k\ge 1$.

\smallskip

{\rm (c)} There is a graphing $G$ such that $\rho_{G_n,r}\to
\rho_{G,r}$ for every $k\ge 1$.
\end{theorem}

Note that in (b) we don't claim uniqueness. We could replace
``graphing'' by ``measure preserving graph''.

A big difference from the dense case is that there does not seem to
be any easy way to construct a sequence that converges to a given
graphing in this sense.

\begin{conj}[Aldous--Lyons]\label{ALD-LYO}
Every graphing is the limit of a convergent sequence of
bounded-degree graphs. Equivalently, every unimodular distribution on
rooted countable graphs with bounded degree is the limit of a bounded
degree graph sequence.
\end{conj}

\subsubsection{Is the limit informative enough?}

The problem of the Regularity Lemma is related to conjecture
\ref{ALD-LYO}. Indeed, suppose that we have a constructive way of
finding, for an arbitrarily large graph $G$ with bounded degree, a
graph $H$ of size bounded by a function of $r$ and $\eps$ that
approximates the distribution of $r$-neighborhoods in $G$ with error
$\eps$. The same construction should also work with a graphing
instead of $G$. Letting $r\to\infty$ and  $\eps\to 0$, this would
give a sequence of finite bounded degree graphs converging to the
given graphing.

Part of the problem is to recognize ``globally'' when $H$ is a good
approximation of $G$. Is there a good notion of ``distance''
(analogous to $\delta_\square$) for graphs with bounded degree?

The limit graphon of a dense sequence of graphs contains very much
information about the asymptotic properties of the sequence. This is
not so for the dense case, unfortunately.

\begin{prob}\label{QU:STR-CONV}
Is there a notion of convergence for graphs with bounded degree that
is stronger than Benjamini--Schramm? (For example, one should be able
to read off from the limit that the graphs are expanders.)
\end{prob}

Let us illustrate this by a couple of simple examples.

\begin{example}\label{EXA:EXPAND}
Let $(G_n)$ be a sequence of 3-regular bipartite expander graphs with
their girth tending to infinity. Let $H_i$ consist of two disjoint
copies of $G_i$. The Benjamini--Schramm limit of both sequences is a
distribution concentrated on a single 3-regular rooted tree. In the
Elek description, we get a graphing $(\Omega, T_1,T_2,T_3)$, where
$T_1,T_2$ and $T_3$ generate a free group which acts on $\Omega$
without fixed points.

This limit graphing is not uniquely determined. One feels that in the
case of the limit of the sequence $(G_n)$, the action of the free
group should ergodic, while in the case of the $H_n$, $\Omega$ should
split into two invariant subsets of measure $1/2$. So it appears that
in the limit object, the underlying $\sigma$-algebra also carries
combinatorial information. This is in stark contrast with the dense
case \cite{BCL}.
\end{example}

\begin{example}\label{EXA:GRIDS}
Let $G_n$ denote the $n\times n$ grid. The Benjamini-Schramm limit
object is a probability distribution concentrated on the infinite
grid with a specified root (the ``origin''). A limit graphing can be
described as the uniform measure on the 2-dimensional torus, together
with the rotations by an irrational number $\alpha$ in one coordinate
and the other.

However, in many respects the ``right'' limit object of the sequence
of grids is a solid square. In other words, instead of larger and
larger pieces of the infinite grid, we consider finer and finer
subdivisions of the unit square.
\end{example}

This last example suggests that we can consider our graphs ``on a
different scale'', and study them as metric spaces with the usual
graph distance as metric, normalized by the diameter. We can then
consider the limit of these metric spaces in the sense of Gromov
\cite{Gro}. For example, the limit of a sequence of larger and larger
square grids in this sense is a (full) square. This global structure
is not revealed by the Benjamini--Schramm limit.

It is easy to construct examples where the interesting structure of
the graphs appears on an intermediate scale. It would be very
interesting to describe and possibly unify limit objects belonging to
different scales. Perhaps we can we understand different limit
objects using ultraproducts, similarly to the work of Elek and
Szegedy in the dense case.

\subsection{Convergence from the right}\label{SP-CONV-R}

While the description of convergent sequences in the bounded degree
case lacks some of the key results that hold in the dense case, most
notably a good notion of distance, we can formulate a result (Borgs,
Chayes, Kahn and Lov\'asz \cite{BCKL}) which shows that convergence
defined in terms of homomorphisms from the left and homomorphisms to
the right are equivalent under some circumstances.

To state this, let us define for every simple graph $G$ and weighted
graph $H$ the quantity
\[
u(G,H)=\frac{\log \hom(G,H)}{|V(G)|},
\]
To see the meaning of $u(G,H)$, consider the case when $H$ is simple.
Then $\hom(G,H)\le q^{|V(G)|}$, and so after taking the logarithm and
dividing by $|V(G)|$, we get a number less than $q$. So $u(G,K_q)$
expresses the freedom (entropy) we have in choosing the image of a
node $v\in V(G)$ in a homomorphism $G\to H$.

For a weighted graph $H$, we define and
\[
\beta_{\max} = \max_{i,j} \beta_{ij},\qquad D(H)=\sum_{i,j\in
V(H)}\frac{\alpha_i\alpha_j}{\alpha_H^2}
\Bigl(1-\frac{\beta_{ij}}{\beta_{\max}}\Bigr).
\]

\begin{theorem}\label{THM:HOM-CONV}
Let $(G_n)$ be a sequence of graphs with maximum degree at most $d$.

\smallskip

{\rm (a)} If $(G_n)$ is convergent, then for every weighted graph $H$
be a weighted graph with $D(H) \le \frac{1}{2d}$, the sequence
$u(G_n,H)$ is convergent.

\smallskip

{\rm (b)} Assume that for every $q\ge 1$ there is an $\eps_q>0$ such
that for every weighted graph $H$ on $q$ nodes with $D(H) \le \eps_q$
the sequence $q(G_n,H)$ is convergent. Then the sequence $(G_n)$ is
convergent.
\end{theorem}

In the special case $H=K_q$ is the complete graph on $q$ nodes
(without loops), we have $D(K_q)=1/q$, and $\hom(G,K_q)$ is the
number of $q$-colorings of $G$. So it follows that if $(G_n)$ is
convergent and $q\ge 2d$, then the number of $q$-colorations grows as
$c^{|V(G_n)|}$ for some $c$. It is easy to see that some condition on
$q$ is needed: for example, if $G_n$ is the $n$-cycle and $q=2$, then
$q(G_n,K_2)$ oscillates between $-\infty$ and $\approx 0$ as a
function of $n$.

\section{Testing}\label{TEST}

What can we learn about a huge graph $G$ from sampling? There are two
related, but slightly different ways of asking this question, {\it
property testing} and {\it parameter estimation}.

\subsection{Sample concentration}

Before discussing these tasks, let us address the following concern:
if we take a bounded size sample from a graph, we can see very
different graphs. For a random graph, for example, we can see
anything. The natural way to use the sample $G[S]$ is to compute some
graph parameter $f(G[S])$. But this parameter can vary wildly with
the choice of the sample, so what information do we get?

The following two theorems assert that every reasonably smooth
parameter of a sample is highly concentrated. (Note: we don't say
anything here about the value of the parameter on the whole graph.)

The first version applies to parameters where smoothness is defined
locally. The proof depends on the theory of martingales (Azuma's
Inequality).

\begin{theorem}\label{CONC-SAMPLE1}
Let $f$ be a graph parameter and assume that $|f(G)-f(G')|\le1 $ for
any two graphs on the same node set which differ only in edges
incident with a single node. Then for every graph $G$ and $1\le k\le
|V(G)|$ there is a value $f_0$ such that if $S\subseteq V(G)$ is a
random $k$-subset, then for every $t>0$,
\[
|f(G[S])-f_0|< \sqrt{2tk}
\]
with probability at least $1-e^{-t}$.
\end{theorem}

The second version applies to parameters which are smooth with
respect to our global distance function. The proof follows from a
modification of the proof of Theorem \ref{LEM:SAMPLE2}.

\begin{theorem}\label{CONC-SAMPLE2}
Let $f$ be a graph parameter and assume that $|f(G)-f(G')|\le
d_\square(G,G')$ for any two graphs on the same node set. Then for
every graph $G$ and $1\le k\le |V(G)|$ there is a value $f_0$ such
that if $S\subseteq V(G)$ is a random $k$-subset, then
\[
|f(G[S])-f_0|<\frac{20}{\sqrt{k}}
\]
with probability at least $1-2^{-k}$.
\end{theorem}

\subsection{Parameter estimation}\label{PARTEST}

We want to determine some parameter of a very large graph $G$. Of
course, we'll not be able to determine the exact value of this
parameter; the best we can hope for is that if we take a sufficiently
large sample, we can find the approximate value of the parameter with
large probability.

To be precise, a graph parameter $f$ is {\it testable}, if for every
$\eps>0$ there is a positive integer $k$ such that if $G$ is a graph
with at least $k$ nodes and we select a set $X$ of $k$ independent
uniform random nodes of $G$, then from the subgraph $G[X]$ induced by
them we can compute an estimate $g(G[X])$ of $f$ such that
\[
\Prob(|f(G)-g(G[X])|>\eps)<\eps.
\]
It is an easy observation that we can always use $g(G[X])=f(G[X])$
(cf. \cite{GT}).

It is easy to see that testability is equivalent to saying that for
every convergent graph sequence $(G_n)$, the sequence of numbers
$(f(G_n))$ is convergent. (So graph parameters of the form $t(F,.)$
are testable by the definition of convergence.) This is, however,
more-or-less just a reformulation of the definition. Paper
\cite{BCLSV1} contains a number of more useful conditions
characterizing testability of a graph parameter. We formulate one,
which is perhaps easiest to verify:

\begin{theorem}\label{THM:PAR-TESTABLE}
A graph parameter $f$ is testable if and only if the following three
conditions hold:

\smallskip

{\rm (i)} For every $\eps>0$ there is an $\eps'>0$ such that if $G$
and $G'$ are two simple graphs on the same node set and
$d_\square(G,G')\le\eps'$ then $|f(G)-f(G')|\le\eps$.

\smallskip

{\rm (ii)} For every simple graph $G$, $f(G(m))$ has a limit as
$m\to\infty$. (Recall that $G(m)$ denotes the graph obtained from $G$
by blowing up each node into $m$ twins. )

\smallskip

{\rm (iii)} If $G^+$ is obtained from $G$ by adding a single isolated
node, then $f(G^+)-f(G)\to 0$ if $|V(G)|\to\infty$.
\end{theorem}

Note that all three conditions are special cases of the statement
that

\smallskip

(iv) {\it if $|V(G_n)|,|V(G'_n)|\to\infty$ and
$\delta_\square(G_n,G'_n)\to 0$, then $f(G_n)-f(G'_n)\to 0$}.

\smallskip

This condition is also necessary, so it is equivalent to its own
three special cases (i)--(iii) in the Theorem.

\begin{example}\label{EXA:MAXCUT}
As a basic example, consider the density of maximum cuts (recall
Section \ref{MULTIWAY}). One of the first substantial results on
property testing \cite{GGR,AKK} is that this parameter is testable.
It is relatively easy to see (using high concentration results like
Azuma's inequality) that if $S$ is a sufficiently large random subset
of nodes of $G$, then $\maxcut(G[S])\ge \maxcut(G)-\eps$: a large cut
in $G$, when restricted to $S$, gives a large cut in $G[S]$. It is
harder, and in fact quite surprising, that if most subgraphs $G[S]$
have a large cut, then so does $G$. This follows from Theorem
\ref{THM:PAR-TESTABLE} above, since conditions (i)--(iii) are easily
verified for $f=\maxcut$.
\end{example}

\begin{example}\label{EXA:FREE}
The free energy \eqref{EQ:FREE-EN} for a fixed weighted graph $H$ is
a more complicated example of a testable parameter, which illustrates
the power of Theorem \ref{THM:PAR-TESTABLE}. It is difficult to
verify directly either the definition, or say condition (iv). The
theorem splits this into three: condition (i) is easy by the
definition of $d_\square(G,G')$; (ii) is a matter of classical
combinatorics, counting mappings that split the twin classes in given
proportions; finally, (iii) is trivial.
\end{example}

\subsection{Dense property testing}\label{PROPTEST-D}

Instead of estimating a numerical parameter, we may want to determine
some property of $G$: Is $G$ 3-colorable? Is it connected? Does it
have a triangle? The answer will of course have some uncertainty. A
precise definition was given by Rubinfeld and Sudan \cite{RS} and
Goldreich, Goldwasser and Ron \cite{GGR}. In the slightly different
context of ``additive approximation'', closely related problems were
studied by Arora, Karger and Karpinski \cite{AKK} (see
e.g.~\cite{Fisch} for a survey). Many extensions deal with situations
where we are allowed to sample more than a constant number of nodes
of the large graph $G$; our concern will be the original setup, where
the sample size is bounded.

A graph property $\PP$ is {\it testable}, if there exists another
property $\PP'$ (called a ``test property'') such that

\smallskip

(a) if a graph $G$ has property $\PP$, then for all $1\le k\le
|V(G)|$ at least $2/3$ of its $k$-node induced subgraphs have
property $\PP'$, and

\smallskip

(b) for every $\eps>0$ there is a $k_\eps\ge 1$ such that if $G$ is a
graph whose edit distance from $\PP$ is at least $\eps |V(G)|^2$,
then for all $k_\eps\le k\le |V(G)|$ at most a fraction of $1/3$ of
the $k$-node induced subgraphs of $G$ have property $\PP'$.

\smallskip

This notion of testability is usually called {\it oblivious testing},
which refers to the fact that no information about the size of $G$ is
assumed. The constants $1/3$ and $2/3$ are arbitrary, and it would
not change the notion of testability if we replaced them by any two
real numbers $0<a<b<1$.

It is surprising that this rather restrictive definition allows many
testable graph properties: for example, bipartiteness,
triangle-freeness, every property definable by a first order formula
\cite{AFKS}.

A surprisingly general result was proved by Alon and Shapira
\cite{ASh}. A graph property $\PP$ is called {\it hereditary} if
$G\in\PP$ implies that $G'\in\PP$ for every induced subgraph $G'$ of
$G$.

\begin{theorem}[Alon--Shapira]\label{THM:HEREDTEST}
Every hereditary graph property is testable.
\end{theorem}

Fischer and Newman \cite{FN} proved that {\it a property is testable
if and only if the normalized edit distance from the property a
testable parameter.} Alon at al.\ characterized testable graph
properties in terms of Szemer\'edi partitions \cite{AFNS}.

Going to the limit gives a tool of studying testability in a
``cleaner'' form (Lov\'asz and Szegedy \cite{LSz4}). It turns out
that this leads to an interesting interplay between the cut-norm and
the $L_1$-norm on $\WW_0$.

A graph property $\PP$ can be thought of as a subset of $\WW_0$
(through the correspondence $G\mapsto W_G$), and we can consider its
closure $\overline\PP$ in the metric space $(\WW_0,\PP)$. For
example, the closure of the set of triangle-free graphs is the set of
{\it triangle-free graphons}, which can be characterized by the
property $t(K_3,W)=0$. More generally, let $\PP$ be a hereditary
graph property. Then its closure is characterized by the (infinitely
many) equations
\begin{equation}\label{EQ:HERED-CLOS}
t_\ind(F,W)=0\qquad\text{for all}\quad F\notin\PP.
\end{equation}

Closures of testable graph properties will be called {\it testable
graphon properties}. These graphon properties can also be
characterized in terms of a sampling method: we consider the
$W$-random graph $\Ge(k,W)$ as the sample of size $k$ from $W$.

\begin{theorem}\label{THM:GRAPHONTEST}
A graphon property $\RR$ is testable if and only if there is a graph
property $\RR'$ such that

\smallskip

(a) $\Pr(\Ge(k,W)\in\RR') \ge 2/3$ for every function $W\in\RR$ and
every $k\ge 1$, and

\smallskip

(b) for every $\eps>0$ there is a $k_\eps\ge 1$ such that
$\Pr(\Ge(k,W)\in\RR') \le 1/3$ for every $k\ge k_\eps$ and every
function $W\in\WW_0$ with $d_1(W,\RR)\ge \eps$.
\end{theorem}

We quote an analytic characterization of testable graphon properties
\cite{LSz4}. Recall that the distances $d_1$ and $d_\square$ are
related trivially by $d_\square\le d_1$. Testability of a property
concerns an inverse relation:

\begin{theorem}\label{PROPTEST-FN}
A graphon property $\RR$ is testable if and only if either one of the
following conditions hold:

\smallskip

{\rm (a)} For every $\eps>0$ there is an $\eps'>0$ such that if
$d_\square(W,\RR)\le \eps'$ for some graphon $W$, then $d_1(W,\RR)\le
\eps$.

\smallskip

{\rm (b)} $d_1(W,\RR)$ is a continuous function of $W$ in the cut
norm.
\end{theorem}

Condition (b) can be viewed as the graphon analogue of the theorem of
Fischer and Newman mentioned above (and the finite theorem can be
derived from it). Condition (a) is a special case of (b).

\begin{example}\label{TEST-FN-EX0}
Let $\RR=\{U\}$, where $U\in\WW$ is the identically $1/2$ function.
Clearly this property is invariant under weak isomorphism. Consider
the random graphs $G_n=\Ge(n,1/2)$; then $\|W_{G_n}-U\|_\square\to 0$
with probability $1$, but $\|W_{G_n}-U\|_1=1/2$ for every $n$. So
this property is not testable by Theorem \ref{PROPTEST-FN}.
\end{example}

Let us sketch how the graphon version of Theorem \ref{THM:HEREDTEST}
follows from this. A property $\RR$ of functions $W\in\WW_0$ is
called {\it flexible} if for every function $U$ such that
$U(x,y)=W(x,y)$ for all $x,y$ with $W(x,y)\in\{0,1\}$, we also have
$U\in\RR$. First, one proves that

\begin{lemma}\label{LEM:CLOS-FLEX}
The closure of a hereditary property is flexible.
\end{lemma}

Indeed, each of the equations \eqref{EQ:HERED-CLOS} is preserved if
we change the value of $W$ at points where this value is positive.

Next, we assume that $\RR$ is a closed flexible property which is not
testable. By Theorem \ref{PROPTEST-FN}, there is a sequence of
functions $W_n$ such that $d_{\square}(W_n,\RR)\to 0$ but
$d_1(W_n,\RR)\geq\eps$ for some fixed $\eps>0$. By Theorem
\ref{THM:COMPACT}, we may assume that $W_n$ converges to some
$W\in\RR$ in the $\|.\|_{\square}$ norm. Let $S_0=W^{-1}(0)$ ,
$S_1=W^{-1}(1)$ and let $Z_n\in\WW_0$ denote the function which is
$1$ on $S_1$, $0$ on $S_0$ and is identical with $W_n$ anywhere else.
By flexibility, we have $Z_n\in\RR$, and by \eqref{EQ:WEAK-ST2},
\[
\|W_n-Z_n\|_1= \int_{S_0}W_n+\int_{S_1}(1-W_n) \to
\int_{S_0}W+\int_{S_1}(1-W)=0  \qquad(n\to\infty),
\]
and so $d_1(W_n,\RR)\to 0$, a contradiction. So it follows that the
closure of every hereditary property is testable.

From this, one can derive that hereditary properties are testable.
There is some further arguments needed, since a graph property can
have a testable closure without itself being testable. (An example is
the property that the graph is complete if the number of nodes is
even but edgeless if the number of nodes is odd.) One can add further
conditions that lead to a characterization, but we don't go into
these technical issues here.

\subsection{Sparse property testing}\label{PROPTEST-S}

We say that a graph property $\PP$ is testable for graphs in $\GG_d$
if for every $\eps>0$ there are integers $r=r(d,\eps)\ge 1$ and
$k=k(d,\eps)$ such that sampling $k$ neighborhoods of radius $r$ from
a graph $G$ with degree bounded by $d$, we can compute ``YES'' or
``NO'' so that:

\smallskip

(a) if we answer ``NO'', then $G\notin\PP$;

\smallskip

(b) if we answer ``YES'', then we can change at most $\eps |V(G)|$
edges in $G$ to get a graph in $\PP$.

An important analogue of the result of Alon and Shapira discussed
above is the following theorem of Benjamini, Schramm and Shapira
\cite{BSchSh}. We must recall a fundamental notion from graph theory:
a {\it minor} of a graph $G$ is any other graph obtained from $G$ by
deleting edges and/or nodes, and contracting edges. A graph property
is {\it minor-closed}, if it is preserved by these operations.
Planarity of a graph is an example of a minor-closed property.

\begin{theorem}\label{THM:BSS}
Every minor-closed property is testable for graphs with bounded
degrees.
\end{theorem}

A related result was proved by Elek \cite{Elek2}:

\begin{theorem}\label{THM:EL-MON}
If a graph property is preserved by edge/node deletion and disjoint
union, then it is testable for graphs with bounded degrees and
subexponential growth.
\end{theorem}

\section{Extremal graph theory}\label{EXT-GRAPHONS}

\subsection{Some classical results}\label{EXT-CLASSIC}

In this section we describe applications of the theory of graph
homomorphisms and graph limits to extremal graph theory. As an
introduction, let us recall some classical results.

Define the {\it Tur\'an graph} $T(n,r)$ ($1\le r\le n$) as follows:
we partition $[n]$ into $r$ classes as equitably as possible, and
connect two nodes if and only if they belong to different classes.

\begin{theorem}[Tur\'an's Theorem]\label{THM:TURAN}
Among all graphs on $n$ nodes containing no $K_k$, the graph
$T(n,k-1)$ has the maximum number of edges.
\end{theorem}

Since we are interested in large $n$ and fixed $k$, the complication
that the classes cannot be exactly equal in size (which causes the
formula for the number of edges of $T(n,k-1)$ to be a bit ugly)
should not worry us. We will be interested in the following
corollary:

\begin{corollary}\label{COR:TURAN}
If a graph on $n$ nodes has more than $\binom{k-1}{2}
\bigl(\frac{n}{k-1}\bigr)^2$ edges, then it contains a $K_k$.
\end{corollary}

The case $k=3$ was proved by Mantel before Tur\'an. We will use this
case to illustrate the ideas, but the general case could be treated
similarly.

One can ask for not just the existence of complete $k$-graphs, but
for their number. Generalizing Tur\'an's Theorem, the following lower
bound was proved by Goodman (for $k=3$) and by Moon and Moser.

\begin{theorem}\label{THM:MOON-MOSER}
If a graph on $n$ nodes has $a\binom{n}{2}$ edges $(0\le a\le 1)$,
then it contains at least $a(2a-1)\dots ((k-2)a-k+1)\binom{n}{k}$
complete $k$-graphs.
\end{theorem}

\begin{figure}[htb]
  \centering
  \includegraphics*[width=.5\textwidth]{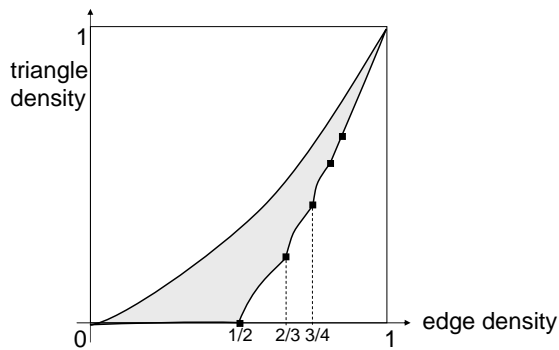}
  \caption{Possible edge and triangle densities of a graph}
  \label{FIG:Turan-plot}
\end{figure}

This bound is tight for Tur\'an graphs, but their edge density
attains only certain values of $a$. The best lower bound in terms of
$a$ and $n$ is quite complicated. To illustrate these complications,
we represent each graph $G$ by the points $(t(K_2,G),t(K_3,G)$ in the
unit square (see Figure \ref{FIG:Turan-plot}). The lower bounding
curve consists of infinitely many concave cubic arcs, and its
validity was only recently proved by Razborov \cite{Razb1}. This was
extended to the best lower bound on the number of $K_4$-s by
Nikiforov \cite{Nik}, but even the edge--$K_q$ diagram is only
conjectural \cite{LSi2} for $q\ge 5$.

One can also ask for an upper bound on the number of complete
$k$-graphs in a graph with given number of edges. A special case of
the Kruskal--Katona Theorem answers this (the whole theorem gives the
precise value, not just asymptotics, and concerns uniform
hypergraphs, not just graphs).

\begin{theorem}\label{THM:KRUS-KAT}
If a graph on $n$ nodes has $a\binom{n}{2}$ edges $(0\le a\le 1)$,
then it contains at most $a^{k/2}\binom{n}{k}$ complete $k$-graphs.
\end{theorem}

Asymptotic equality is attained when the graph consists of a clique
and isolated nodes. Not every edge density $a$ can be realized by
such graphs, but the attainable edge densities are dense in $[0,1]$,
and so Theorem \ref{THM:KRUS-KAT} is asymptotically tight for all
values of $a$.

Instead of counting complete graphs, we one can consider the number
of copies of some other graph $F$ in $G$. We have already come across
counting $4$-cycles twice: in Section \ref{QUASIRAND} and in Section
\ref{GR-OPT}. Giving just the simpler asymptotic version:

\begin{theorem}[Erd\H{o}s]\label{THM:4CYC}
If a graph on $n$ nodes has $a\binom{n}{2}$ edges $(0<a\le 1)$, then
it contains at least $(\frac18+o(1))a^4 n^4$ $4$-cycles.
\end{theorem}

Graphs with asymptotic equality here are quasirandom graphs.

The number of paths of length $k$ is a more difficult question, but
it turns out to be equivalent to a theorem of Blakley and Roy
\cite{BR} in matrix theory. Again asymptotically,

\begin{theorem}\label{THM:BL-ROY}
If a graph on $n$ nodes has $a\binom{n}{2}$ edges $(0<a\le 1)$, then
it contains at least $(\frac12 +o(1)) a^{k-1} n^k$ paths of length
$k$.
\end{theorem}

Regular graphs give asymptotic equality here.

\subsection{Algebraic proofs of extremal graph
results}\label{ALG-PROOFS}

The classical extremal problems in the previous section can be
expressed as algebraic inequalities between the subgraph densities
$t(F,W)$ that hold for all graphons $W$. Often ``going to the
infinity'' provides cleaner formulations (no error terms). Here are a
few examples:

\begin{example}\label{EXA:EXTREM}\strut

(a) {\it Tur\'an's theorem.} We state just the case of triangles (due
to Mantel):
\begin{equation}\label{EQ:MANTEL}
t(K_3,W)=0 \Rightarrow t(K_2,W)\le 1/2,
\end{equation}
which follows from the algebraic inequality due to Goodman
\cite{Good}:
\begin{equation}\label{EQ:GOODMAN}
t(K_3,W)\ge t(K_2,W)(2t(K_2,W)-1).
\end{equation}

\smallskip

(b) {\it The Kruskal--Katona theorem for graphs}:
\begin{equation}\label{EQ:KR-KAT}
t(K_3,W)\le t(K_2,W)^{3/2}.
\end{equation}

\smallskip

(c) {\it Erd\H{o}s's bound on the number of quadrilaterals}:
\begin{equation}\label{EQ:ERD-QUAD}
t(C_4,W)\ge t(K_2,W)^4.
\end{equation}

\smallskip

(d) {\it The Blakley--Roy inequality}:
\begin{equation}\label{EQ:BL-ROY}
t(P_k,W)\ge t(K_2,W)^{k-1}.
\end{equation}

\smallskip

(e) {\it The Sidorenko Conjecture} (unsolved) generalizes the last
two results in the direction that for every bipartite graph $F$,
\begin{equation}\label{EQ:SIDOR}
t(F,W) \ge t(K_2,W)^{|E(F)|}.
\end{equation}
This conjecture is proved for trees, many small graphs, complete
bipartite graphs (Sidorenko \cite{Sid}) and also for cubes (Hatami
\cite{Hat}).
\end{example}

Using the formalism introduced above, the results in example
\ref{EXA:EXTREM} can be expressed as follows:
\begin{align*}
{\rm (a)}&\qquad K_3 \ge 2 K_2{}^2 - K_2;\\
{\rm (b)}&\qquad K_2{}^3 \ge K_3{}^2 ;\\
{\rm (c)}&\qquad C_4 \ge K_2{}^4;\\
{\rm (c)}&\qquad P_4 \ge K_2{}^3 ;\\
{\rm (d)}&\qquad F \ge K_2{}^{|E(F)|} \qquad\text{(if $F$ is
bipartite)}.
\end{align*}

The first three inequalities can be proved easily using the
reflection positivity of the graph parameters $t(.,W)$. We will
illustrate the method by deriving (a) through formal algebraic
manipulations.

\begin{figure}[htb]
  \centering
  \bigskip
  \includegraphics*[width=.7\textwidth]{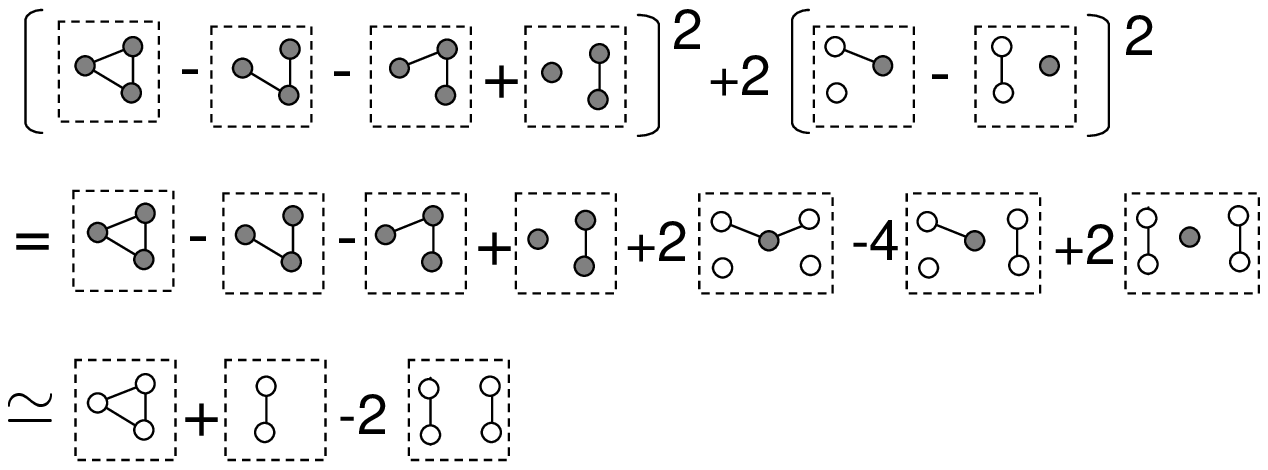}
  \caption{A computation proving the Mantel-Tur\'an Theorem}
  \label{FIG:TURAN}
\end{figure}

\medskip\noindent{\bf Proof of (a)} (Goodman's extension of the
Mantel--Tur\'an Theorem). Let $F$ denote the graph $K_2K_1$ (an edge
and an isolated node), and let $F_1$, $F_2$ and $F_3$ be obtained
from $F$ by labeling all three nodes, one endpoint of the edge, and
the isolated node, respectively. Consider the quantum graph
$\widetilde{F_1}{}^2 + 2 (F_2-F_3)^2$, which is obviously
nonnegative. Unlabeling the nodes and deleting isolated nodes, we get
$K_3-2K_2{}^2+K_2$, which is thus nonnegative (see Figure
\ref{FIG:TURAN}).

Of the above inequalities, also (b) and (c) can be proved by similar
arguments. The Blakley-Roy inequality (c) is more difficult, but some
extension of this kind of argument does work \cite{Kun}. Sidorenko's
conjecture (d) would of course be very nice to prove this way (or by
any other means).

Using related methods, Razborov \cite{Razb1} solved the long-standing
problem of characterizing the possible (edge-density,
triangle-density) pairs, which in this setting means a description of
the set $(t(K_2,W),t(K_3,W)):~W\in\WW_0)$ by algebraic inequalities.

The inequality in (c) also follows from reflection positivity if $k$
is even. It is not known whether (c) for odd $k$ (or perhaps every
valid algebraic inequality between subgraph densities) follows from a
finite number of semidefiniteness inequalities. However, every valid
linear inequality between homomorphism densities follows from
semidefiniteness constraints (equivalently, from ``sums of squares''
computations in graph algebras), as we shall see in the next section.

\subsection{Positivstellensatz for graphs and spectral
norms}\label{POSITIVST}

The machinery introduced in the previous sections allows us to
suggest a very general approach to extremal graph theory.

We can define the following partial order on $\QQ_0$: we say that a
quantum graph $x\ge 0$, if $t(x,W)\ge 0$ for all $W\in\WW_0$.

Let us call a quantum graph $y$ a {\it square-sum} if there are
$k$-labeled quantum graphs $y_1,\dots,y_k$ for some $k$ such that $y$
can be obtained from $\sum_i y_i^2$ by forgetting the labels. It is
easy to see that every square-sum satisfies $y\ge 0$.

As an example, recall the definition \eqref{EQ:SZITA} of the
``inclusion-exclusion'' quantum graph $\widehat{F}$. Let us label all
nodes of $\widehat{F}$, square it, and then forget the labels: we
obtain $\widehat{F}$ itself. This implies that $\widehat{F}\ge 0$ for
all $F$. In the special case when $W=W_G$ for some graph $G$, this
also follows from our previous remark that $t(\widetilde{F},G)$ is a
probability, and hence nonnegative.

Is there a quantum graph $x\ge 0$ which is not a square sum? I
suspect that such quantum graphs exist, but it might be difficult to
prove this property. However, the following weaker result can be
proved \cite{LSz7}.

\begin{theorem}\label{THM:W-POS}
Let $x$ be a quantum graph. Then $x\ge 0$ if and only if for every
$\eps>0$ there is a square-sum $y$ such that $N(y)\le N(x)$ and
$\|x-y\|_2<\eps$.
\end{theorem}

The proof depends on the duality theory of semidefinite programs.
Note that we do not claim that the $k$-labeled quantum graphs $y_i$
in the square-sum representation of $y$ also have bounded $N(y_i)$;
the proof gives arbitrarily large graphs if $\eps$ is small.

In analogy with the Positivstellensatz for real polynomials, we may
try to represent quantum graphs $x\ge0$ as quotients of square-sums:
if $y$ and $z$ are square-sums and $y=zx+x$, then $x\ge 0$.

We mention a couple of related questions. For every even positive
integer $k$, the functional $t(C_k,W)^{1/k}$ defines a norm on $\WW$
(the Neumann-Schatten norm). This suggests the question: For which
other simple graphs $F$ is $t(F,W)^{1/|E(F)|}$ a norm (or seminorm)
on $\WW$? Hatami \cite{Hat} proved that if a simple graph $F$ has the
property that $\|W\|=t(F,|W|)^{1/|E(F)|}$ is a norm, then it
satisfies Sidorenko's conjecture \ref{EXA:EXTREM}(d). He also proved
that all cubes have this property.

In view of the usefulness of extending graphs to graphons, it seems
natural to define graph algebras of infinite linear combinations of
graphs with appropriate convergence properties. It is not worked out,
however, what the structure of the resulting algebra is, and how it
is related to graphons.

\subsection{The maximum distance from a hereditary graph property}

A surprisingly general result is the theorem of Alon and Stav
\cite{ASt}, proving that for every hereditary property, a random
graph with appropriate density is asymptotically the farthest from
the property in edit distance. The analytic results developed in this
paper allow us to state and prove a simple analytic analogue of this
fact, from which the original result follows along with
generalizations.

\begin{theorem}[Alon and Stav]\label{ALON-STAV}
For every hereditary graph property $\PP$ there is a number $p$,
$0\le p\le 1$, such that for every graph $G$ with $|V(G)|=n$,
\[
d_1(G,\PP)\le \E(d_1(\Ge(n,p),\PP)) + o(1) \qquad (n\to\infty).
\]
\end{theorem}

The following theorem \cite{LSz4} states a graphon version of this
fact.

\begin{theorem}\label{ALON-STAV-FN}
If $\RR$ is the closure of a hereditary graph property, then the
maximum of $d_1(.,\RR)$ is attained by a constant function.
\end{theorem}

Our point in giving this generalization is to illustrate the power of
extending graph problems to a continuum. The key observation is the
following, which follows from Lemma \ref{LEM:CLOS-FLEX}.

\begin{lemma}\label{LEM:CONV}
If $\RR$ is the closure of a hereditary graph property, then the set
$\WW_0\setminus \RR$ is convex.
\end{lemma}

Hence it follows that the $d_1$ distance from $\PP$ is a concave
function on $\WW_0\setminus \RR$. Since $\WW_0\setminus \RR$ is
obviously invariant under the group of invertible measure preserving
transformations of $[0,1]$, it is not hard to argue that there is a
point (graphon) in $\WW_0\setminus \RR$ maximizing the distance from
$\PP$ which is invariant under these measure preserving
transformations, and so it must be a constant function.

\subsection{Which graphs are extremal? (Finitely forcible graphons)}

We call a graphon $W\in\WW_0$ {\it finitely forcible} if there exist
a finite list of graphs $F_1,\dots, F_m$ and real numbers
$a_1,\dots,a_m$ such that the equations $t(F_1,U)= a_1,\dots,t(F_m,U)
=a_m$ are satisfied by precisely those functions $U\in\WW_0$ which
arise from $W$ by measure preserving transformations.

Let us consider a very general type of graph theoretic extremal
problem:
\begin{align}\label{EQ:EXT}
\text{~maximize}~&t(f,W)\nonumber\\
\text{subject to}~&t(g_1,W)=a_1\nonumber\\
&t(g_2,W)=a_1\\
&~\vdots\nonumber\\
&t(g_k,W)=a_1\nonumber
\end{align}
where $f,g_1,\dots,g_k$ are given quantum graphs. Most of the graphon
versions of extremal problems discussed so far fit in this scheme.

It is easy to see that every finitely forcible graphon is the
solution of an extremal problem of the type \eqref{EQ:EXT}. We
conjecture the following converse:

\begin{conj}\label{PROB:FINFORCE-EXT}
Every extremal problem has a finitely forcible optimum. In other
words, if a finite set of constraints of the form $t(F_i,W)=a_i$ is
satisfied by some graphon, then it is satisfied by a finitely
forcible graphon.
\end{conj}

This may seem far fetched, but the following heuristic supports it.
Suppose that $t(F_1,W)=a_1,\dots,t(F_k,W)=a_k$ has a solution in $W$,
but this is not forced by these constraints. Then there is a graph
$F$ such that $t(F,W)$ is not determined, i.e., $a=\min t(F,W)<\max
t(F,W)=b$ (the max and min are taken over all solutions $W$ of the
system). Now add one of the conditions $t(F,W)=a$ or $t(F,W)=b$ to
the system and repeat. It seems that in very few (2-3) steps we
always get a unique solution, i.e., a finitely forcible graphon.

Almost all classical extremal problems have a solution that is a
stepfunction. It was shown by Lov\'asz and S\'os \cite{LSos} that
{\it every stepfunction is finitely forcible}, and it was conjectured
that these are the only ones. Recently B.~Szegedy and Lov\'asz
\cite{LSz6} found other finitely forcible graphons, and so the
problem of characterizing finitely forcible graphons is wide open.

We mention two examples of finitely forcible graphons that are not
stepfunctions (the proof is not quite easy).

\begin{example}\label{EXA:WPOLY}
Let $p(x,y)$ is a symmetric real polynomial that is monotone
increasing on $[0,1]^2$. Define
\[
W(x,y)=
  \begin{cases}
    1, & \text{if $p(x,y)>0$}, \\
    0, & \text{otherwise},
  \end{cases}
\]
Then $W$ is finitely forcible. It is conjectured that monotonicity is
not needed here.

In contrast, one can show that if $W(x,y)$ is a polynomial in $x$ and
$y$ (not a function of the sign), then it is {\it not} finitely
forcible.
\end{example}

\begin{example}\label{EXA:WFRACT}
Let
\[
W(x,y)=
  \begin{cases}
    1, & \text{if the first bit where the binary expansions of $x$ and $y$ differ}\\
       & \text{is at an odd position}, \\
    0, & \text{otherwise},
  \end{cases}
\]
The $W$ is finitely forcible.
\end{example}

\end{document}